\documentclass[journal]{IEEEtran}      

\pdfoutput=1
\IEEEoverridecommandlockouts                              

\newtheorem{Remark}{Remark}

\newtheorem{Definition}{Definition}

\newenvironment{Proof}{\noindent{\em Proof:\/}}{\hfill $\Box$\par}
\newtheorem{Theorem}{Theorem}
\newtheorem{Lemma}{Lemma}

\newtheorem{Assumption}{Assumption}

\usepackage{graphicx}
\usepackage{amsmath,bm}
\usepackage{amssymb}
\usepackage{algorithm}
\usepackage{algpseudocode}
\usepackage{cite}
\usepackage{hyperref}
\usepackage[caption=false]{subfig}
\hypersetup{hidelinks}
\allowdisplaybreaks
\usepackage{color}

\algdef{SE}[DOWHILE]{Do}{doWhile}{\algorithmicdo}[1]{\algorithmicwhile\ #1}%

\newcommand{\mathactivatecomma}{%
  \begingroup\lccode`~=`\,
  \lowercase{\endgroup\edef~}{\mathchar\the\mathcode`\,\penalty0 }}

\algnewcommand{\Initialize}[1]{%
  \State \textbf{Initialize: $j\in \mathcal{V}^i, i \in \mathcal{I}$}
  \Statex \hspace*{\algorithmicindent}\parbox[t]{.8\linewidth}{\raggedright #1}
}

\algnewcommand{\Iteration}[1]{%
  \State \textbf{Iteration $(t\geq 0)$: $j\in \mathcal{V}^i, i \in \mathcal{I}$}
  \Statex \hspace*{\algorithmicindent}\parbox[t]{.8\linewidth}{\raggedright #1}
}

\algnewcommand{\Output}[1]{%
  \State \textbf{Output: $j\in \mathcal{V}^i, i \in \mathcal{I}$}
  \Statex \hspace*{\algorithmicindent}\parbox[t]{.8\linewidth}{\raggedright #1}
}

\title{\LARGE \bf
Nash Equilibrium Seeking in $N$-Coalition Games via a Gradient-Free Method
}

\author{Yipeng Pang and Guoqiang Hu
\thanks{This research was supported by Singapore Ministry of Education Academic Research Fund Tier 1 RG180/17(2017-T1-002-158).}
\thanks{Y. Pang and G. Hu are with the School of Electrical and Electronic Engineering, Nanyang
Technological University, 639798, Singapore
        {\tt\small ypang005@e.ntu.edu.sg, gqhu@ntu.edu.sg}.}%
}

\begin{document}

\bstctlcite{IEEEexample:BSTcontrol}

\maketitle
\thispagestyle{empty}
\pagestyle{empty}

\begin{abstract}
This paper studies an $N$-coalition non-cooperative game problem, where the players in the same coalition cooperatively minimize the sum of their local cost functions under a directed communication graph, while collectively acting as a virtual player to play a non-cooperative game with other coalitions. Moreover, it is assumed that the players have no access to the explicit functional form but only the function value of their local costs. To solve the problem, a discrete-time gradient-free Nash equilibrium seeking strategy, based on the gradient tracking method, is proposed. Specifically, a gradient estimator is developed locally based on Gaussian smoothing to estimate the partial gradients, and a gradient tracker is constructed locally to trace the average sum of the partial gradients among the players within the coalition. With a sufficiently small constant step-size, we show that all players' actions approximately converge to the Nash equilibrium at a geometric rate under a strongly monotone game mapping condition. Numerical simulations are conducted to verify the effectiveness of the proposed algorithm. 
\end{abstract}


\section{Introduction}
Recently, great efforts have been devoted to the research of collaboration and competition among multiple rational decision-makers, of which distributed optimization and Nash equilibrium (NE) seeking in non-cooperative games are the two main lines. 
In particular, distributed optimization problem is concerned with a network of agents that cooperatively minimize a combination of their local cost functions, which has shown great interests in the applications of parameter estimation, source localization, resource allocation, multi-robot coordination, \textit{etc}. Non-cooperative game, on the other hand, considers a group of players that are self-interest motivated to minimize their own individual cost function in response to other players' actions, which has been widely applied in the fields of transportation network control, power network control, electricity markets, smart grids, \textit{etc}. 

This paper subsumes the research of both distributed optimization and NE seeking in non-cooperative games, by considering an $N$-coalition non-cooperative game. In this game, each coalition consists of a number of players. On the one hand, each coalition can be regarded as a virtual player that aims to minimize its cost function by adjusting its actions based on other coalition's actions in a non-cooperative game. On the other hand, the coalition's cost function is defined as the sum of the local cost functions associated to the players in the corresponding coalition, and minimization of such cost is realized by the collaboration among the corresponding players. Moreover, it is assumed that the players have no access to the explicit functional form but only the function value of their local costs. A discrete-time gradient-free NE seeking strategy is developed to drive the players' actions to the NE.

\textbf{Related Work}. 
Recently, vast results on the design of NE seeking algorithms have been reported to solve matrix games, potential games, aggregate games, generalized games \cite{Grammatico2017,Deng2019a,Yi2019,Pavel2020}, to list a few. 
For the $N$-coalition game problem, the formulation is related to the problem of two sub-networks zero-sum games (\textit{e.g.}, \cite{Gharesifard2013,Lou2016}), where two adversarial networks have opposing objectives with regard to the optimization of a common cost function, being collaboratively minimized by the agents in the corresponding network. Then, the $N$-coalition game problem is essentially an extension of such problem with the number of subnetworks being $N$. 
To solve the problem, the work in \cite{Ye2018} proposed a NE seeking strategy, where a dynamic average consensus protocol was adopted to estimate the averaged gradients of the coalitions' cost functions, and a gradient play method was implemented to drive the states to the equilibrium. It was further improved in \cite{Ye2019} to reduce the communication and computation costs by incorporating an interference graph which characterizes the interactions among the players in each coalition. 
The $N$-coalition game problem was also studied in \cite{Zeng2019}, where the players are subject to both local set constraints and a coupled inequality constraint, and the associated cost functions are allowed to be non-smooth. A distributed algorithm with the use of projected operators, subgradient dynamics, and differential inclusions was proposed to find the generalized NE. Different from all these works, the authors in \cite{Ye2020} proposed an extremum seeking-based approach without the knowledge on the explicit expressions of the players' local cost functions, which is the most relevant to the problem considered in this paper. However, it should be noted that the methods proposed for $N$-coalition game problems in the existing literature are continuous-time approaches. This paper intends to devise a discrete-time gradient-free NE seeking strategy to drive the players' actions to the NE.

The newly proposed gradient-free NE seeking strategy follows the idea of Gaussian smoothing to estimate the gradient of a function, which was firstly reported in \cite{Nesterov2017} in convex optimization, and further studied in both distributed optimization \cite{Yuan2015,Pang2017,Pang2018,Pang2019b,Pang2020} and non-cooperative games \cite{Pang2020a}. This type of methods, including payoff-based methods in \cite{Tatarenko2019,Tatarenko2019a} and extremum seeking-based methods in \cite{Ye2020}, can be regarded as non-model based approaches, since the implementation of such algorithms does not require the model information. One advantage of the Gaussian smoothing compared with other non-model based methods is that it can deal with both smooth and non-smooth problem setups. Though this technique has been well studied in \cite{Nesterov2017} and its extensions in \cite{Yuan2015,Pang2017,Pang2018,Pang2019b,Pang2020} for optimization problems, the derived properties cannot be directly applied to non-cooperative game problems due to the couplings in the players' actions. Moreover, the algorithm based on Gaussian smoothing technique can only drive the players' actions to the NE of a smoothed game but not the original game. How close between the NEs of the smoothed and the original games has not been studied yet.
With the Gaussian smoothing technique, the proposed strategy utilizes the estimated gradient information in a gradient tracking method to trace the gradient of the coalition's cost function. The gradient tracking technique was widely adopted in distributed optimization to achieve a fast rate of convergence with a constant step-size, see \cite{Xi2016a,Xin2019b,Pu2020,Xin2019}, but receives little attention in NE seeking. 
It should be highlighted that solving $N$-coalition game problems is not equivalent to solving $N$ independent distributed optimization problems among $N$ coalitions, since these $N$ distributed optimization problems are coupled with each other in terms of the players' actions, and each coalition has only authority of its own action rather than the full authority of the entire decision variables in distributed optimization problems. The action of one coalition will have influence to the actions of other coalitions. Moreover, the cost functions in non-cooperative games are only convex with respect to the player's own action, rather than the entire decision variable. Such partial convexity also brings in additional challenges in the convergence analysis of the gradient tracking method in non-cooperative games.

\textbf{Contributions}. 
The major contributions of this paper are threefold.
1). Different from the existing works on $N$-coalition game problems \cite{Ye2018,Ye2019,Zeng2019}, the problem considered in this paper follows the settings in \cite{Ye2020}, where the players are assumed to have no access to the explicit functional form, but only the function value of their local cost functions.
As the methods proposed for $N$-coalition game problems in all existing literature \cite{Ye2018,Ye2019,Zeng2019,Ye2020} are continuous-time approaches, this paper is the first attempt to address $N$-coalition game problems with a discrete-time method, where a non-model based NE seeking strategy is proposed via a gradient-free method. 
2). This paper adopts the Gaussian smoothing technique to estimate the partial gradient of the player's local cost function, and a gradient tracking method is employed to trace the gradient of the coalition's cost function. Since the Gaussian smoothing-based algorithms can only achieve convergence to the NE of a smoothed game but not the original game, we are the first to explicitly quantify the distance between the NEs of the smoothed and the original games as a function of the smoothing parameter.
As compared to the existing gradient tracking methods, such as \cite{Xi2016a,Xin2019b,Pu2020,Xin2019}, this work is the first attempt to adopt gradient tracking techniques in $N$-coalition games.
3). The convergence property of the proposed algorithm is carefully studied. Under the standard assumptions on the graph connectivity, local cost functions and game mappings, it is derived that the players' actions driven by the proposed algorithm with a small constant step-size are able to converge to a small neighborhood of the NE at a geometric rate with the error being proportional to the step-size and the smoothing parameter. 

\textbf{Notations}.
We use $\mathbb{R}$, $\mathbb{R}_{\geq0}$ and $\mathbb{R}^p$ to denote real numbers, non-negative real numbers, and $p$-dimensional column vectors, respectively. $\mathbf{1}_p$ ($\mathbf{0}_p$) represent a vector with all elements equal to $1$ ($0$), and $I_p$ denotes the $p\times p$ identity matrix. 
For a vector $u$, we use $\|u\|$ for its standard Euclidean norm, \textit{i.e.}, $\|u\| \triangleq \sqrt{\langle u, u\rangle}$
For a vector $a$ or a matrix $A$, we use $[a]_i$ to denote its $i$-th entry, and $[A]_{ij}$ to denote its element in the $i$-th row and $j$-th column. The transpose and spectral norm of a matrix $A$ are denoted by $A^\top$ and $\|A\|$, respectively. We use $\rho(A)$ to represent the spectral radius of a square matrix $A$. 
For a totally differentiable function $f(x,y)$, we use $\nabla f(x,y)$ for its total derivative, and $\nabla_xf(x,y)$ for its partial derivative with respect to $x$ for a fixed $y$.

\section{Problem Statement}\label{sec:problem_formulation}

\subsection{Game Formulation}


\begin{Definition}\label{definition_N_coalition}
($N$-coalition Games). An $N$-coalition non-cooperative game is defined by $\Gamma(\mathcal{I},\{f^i\},\{\Omega^i\})$, where each coalition $i$ (indexed by $\mathcal{I} \triangleq \{1, 2, \ldots, N\}$) owns a cost function $f^i$, and consists of $n_i$ number of players (denoted by $\mathcal{V}^i \triangleq \{1,2,\ldots,n_i\}$) with each player's action subject to a constrained set $\Omega^i_j\subset\mathbb{R}^{p^i_j}$, where $p^i_j$ denotes the dimension of the action of player $j$ coalition $i$. Denote $n \triangleq \sum_{i=1}^Nn_i$, $p_i\triangleq\sum_{j=1}^{n_i}p^i_j$, $p \triangleq \sum_{i=1}^{N}p_i$, $\Omega^i=\Omega^i_1\times\cdots\times\Omega^i_{n_i}\subset\mathbb{R}^{p_i}$, and $\Omega \triangleq \Omega^1\times\cdots\times\Omega^N\subset\mathbb{R}^p$. The cost function $f^i:\Omega\to\mathbb{R}$ is defined as, 
\begin{equation*}
f^i(\mathbf{x}^i,\mathbf{x}^{-i}) \triangleq \frac1{n_i}\sum_{j=1}^{n_i}f^i_j(\mathbf{x}^i,\mathbf{x}^{-i}), \quad\mathbf{x}^i \in \Omega^i, \quad \forall i\in\mathcal{I},
\end{equation*}
where $f^i_j(\mathbf{x}^i,\mathbf{x}^{-i})$ is a local cost function of player $j$ in coalition $i$, $\mathbf{x}^i \triangleq [x^{i\top}_1,\ldots,x^{i\top}_{n_i}]^\top \in \Omega^i$ is the collection of all players' actions in coalition $i$, and $\mathbf{x}^{-i}\in\Omega\backslash\Omega^i$ is the collection of all players' actions other than coalition $i$. $x^i_j\in\Omega^i_j$ is the action of player $j$ in coalition $i$. Collectively, we denote $\mathbf{x}\triangleq (\mathbf{x}^{i},\mathbf{x}^{-i}) = [\mathbf{x}^{1\top},\ldots,\mathbf{x}^{N\top}]^\top$. 
\end{Definition}

\begin{Definition}
(Nash Equilibrium of $N$-Coalition Games). A vector $\mathbf{x}^* \triangleq (\mathbf{x}^{i*},\mathbf{x}^{-i*})\in\Omega$ is said to be a Nash equilibrium of the $N$-coalition non-cooperative game $\Gamma(\mathcal{I},\{f^i\},\{\Omega^i\})$, if and only if
\begin{align*}
f^i(\mathbf{x}^{i*},\mathbf{x}^{-i*})\leq f^i(\mathbf{x}^i,\mathbf{x}^{-i*}),\quad \forall x^i\in\Omega^i,\quad \forall i\in \mathcal{I}.
\end{align*}
\end{Definition} 


In this paper, we consider $\Omega^i_j = \mathbb{R}$ for simplicity. For an $N$-coalition non-cooperative game $\Gamma(\mathcal{I},\{f^i\},\{\mathbb{R}^{n_i}\})$, the players in the same coalition cooperatively minimize the summation of their local cost functions, while collectively acting as a virtual player to play a non-cooperative game with other coalitions. To achieve the mutual interest inside the coalition, the players in each coalition $i\in\mathcal{I}$ are equipped with a communication network, characterized by a directed graph $\mathcal{G}_i(\mathcal{V}^i,\mathcal{E}^i)$ with an adjacency matrix $A^i\in\mathbb{R}^{n_i\times n_i}$, where $[A^i]_{jk}>0$ if $(k,j)\in\mathcal{E}^i$ and $[A^i]_{jk}=0$ otherwise. We assume $(k,k)\in\mathcal{E}^i, \forall k\in\mathcal{V}^i$.
Moreover, it is assumed that all players have limited knowledge on their local cost functions, similar to the settings in \cite{Pang2020a,Ye2020}. That is, the player can only access the output of the local cost function, whose explicit form is assumed to be unknown.
Suppose that the considered $N$-coalition game $\Gamma(\mathcal{I},\{f^i\},\{\mathbb{R}^{n_i}\})$ admits a NE.
The objective of this paper is to design a NE seeking strategy such that all players' actions converge to a NE.

This paper considers the case where each agent has full access to all other agents' decisions, which is known as \textit{full-decision information}. In this game setup, even though agents have full access to other agents' decisions, the agents' cost functions including the gradient information are still private to other agents. Hence, across the coalitions, the agents do not need to have extra communication. In the same coalition, since the agents need to collaboratively minimize the coalition's cost function, so they need to communicate with each other to obtain some necessary information, such as the gradient tracker in this paper. Thus, we make the following assumptions on the communication graph of each coalition and the players' local cost functions.

\begin{Assumption}\label{assumption_graph}
The digraph $\mathcal{G}_i$ is strongly connected. The associated adjacency matrix $A^i$ is doubly-stochastic, \textit{i.e.}, $\mathbf{1}_{n_i}^\top A^i = \mathbf{1}_{n_i}^\top$ and $A^i\mathbf{1}_{n_i} = \mathbf{1}_{n_i}$.
\end{Assumption}

Define $\sigma_{A^i} \triangleq \|A^i-\frac1{n_i}\mathbf{1}_{n_i}\mathbf{1}_{n_i}^\top\|$. It follows from \cite[Lemma~1]{Pu2020} that $\sigma_{A^i}<1$. Denote by $\bar{\sigma} \triangleq \max_{i\in\mathcal{I}}\sigma_{A^i}$ and $\varsigma \triangleq \max_{i\in\mathcal{I}}(1+\sigma_{A^i}^2)/(1-\sigma_{A^i}^2)$.

\begin{Assumption}\label{assumption_local_f_lipschitz}
For each $j\in\mathcal{V}^i,i\in\mathcal{I}$, the local cost function $f^i_j(\mathbf{x}^i,\mathbf{x}^{-i})$ is convex in $\mathbf{x}^i$, and continuously differentiable in $\mathbf{x}$. The total gradient $\nabla f^i_j(\mathbf{x})$ is $\mathcal{L}$-Lipschitz continuous in $\mathbf{x}$, \textit{i.e.}, $\|\nabla f^i_j(\mathbf{x}) - \nabla f^i_j(\mathbf{x}')\|\leq \mathcal{L}\|\mathbf{x} - \mathbf{x}'\|$, $\forall \mathbf{x}, \mathbf{x}'\in\mathbb{R}^n$.
\end{Assumption}




Next, we define the game mapping of $\Gamma(\mathcal{I},\{f^i\},\{\mathbb{R}^{n_i}\})$ as given by $F(\mathbf{x})\triangleq [\nabla_{\mathbf{x}^1}f^1(\mathbf{x})^\top,\ldots,\nabla_{\mathbf{x}^N}f^N(\mathbf{x})^\top]^\top$.

In this paper, we impose a strong monotonicity assumption on the game mapping $F(\mathbf{x})$, which is common in many works, \textit{e.g.}, \cite{Deng2019a,Salehisadaghiani2016,Tatarenko2018,Pang2020a,Yi2019}.
\begin{Assumption}\label{assumption_game_mapping}
The game mapping $F(\mathbf{x})$ of game $\Gamma$ is $\chi$-strongly monotone, \textit{i.e.}, for $\forall\mathbf{x},\mathbf{x}'\in\mathbb{R}^n$, we have $\langle F(\mathbf{x})-F(\mathbf{x}'), \mathbf{x}-\mathbf{x}' \rangle\geq \chi\|\mathbf{x}-\mathbf{x}'\|^2$.
\end{Assumption}
\begin{Remark}
It is known that under Assumptions~\ref{assumption_local_f_lipschitz} and \ref{assumption_game_mapping}, game $\Gamma$ admits a unique Nash equilibrium \cite[Thm.~3]{Scutari2014}.
\end{Remark}

\subsection{Preliminaries on Gaussian Smoothing}
To facilitate the gradient-free techniques, some preliminaries on Gaussian smoothing and randomized gradient oracle originated from \cite{Nesterov2017} are presented in this part. 

For the local cost function $f^i_j(\mathbf{x})$ of player $j\in\mathcal{V}^i$ in coalition $i\in\mathcal{I}$, a Gaussian-smoothed function $f^i_{j,\mu}(\mathbf{x})$ with respect to $\mathbf{x}$ is defined as
\begin{align}
  f^i_{j,\mu}(\mathbf{x}) \triangleq \frac1{\kappa}\int_{\mathbb{R}^n} f^i_j(\mathbf{x}+\mu\xi)e^{-\frac12\|\xi\|^2}d\xi, \label{eq:smooth_function_definition}
\end{align}
where $\kappa \triangleq \int_{\mathbb{R}^n} e^{-\frac12\|\xi\|^2}d\xi = (2\pi)^{n/2}$, $\xi \in \mathbb{R}^n$ is a normally distributed random variable, and $\mu \geq 0$ is called the smoothing parameter. 

Within each coalition $i\in\mathcal{I}$, the randomized gradient-free oracle for player $j$ with respect to player $k$, $j,k\in\mathcal{V}^i$, at step $t\geq 0$ is designed as\footnote{In some cases, the two-point sampling in Gaussian smoothing may necessitate some coordination between players in the gradient approximation, which may be more stringent than the one-point sampling in other non-model based methods.}
\begin{align}
\pi^i_{jk}(\mathbf{x}_t) \triangleq \frac{f^i_j(\mathbf{x}_t+\mu\xi^i_{j,t})-f^i_j(\mathbf{x}_t)}{\mu}[\xi^i_{j,t}]^i_k, \label{grad_oracle}
\end{align}
where $[\xi^i_{j,t}]^i_k$ denotes the $(\sum_{l=0}^{i}n_l+k)$-th element of $\xi^i_{j,t}$ with $n_0 = 0$ and $\xi^i_{j,t}$ being its own version of $\xi$ at step $t$, and $\mu>0$.
Let $\mathcal{F}_t$ denote the $\sigma$-field generated by the entire history of all random variables from step $0$ to $t-1$.
Then, the following properties on $f^i_{j,\mu}(\mathbf{x})$ and $\pi^i_{jk}(\mathbf{x}_t)$ are directly collected from \cite{Nesterov2017}.

\begin{Lemma}\label{lemma:property_f_mu}
(see \cite{Nesterov2017}) Under Assumption~\ref{assumption_local_f_lipschitz}, for $\forall j,k\in\mathcal{V}^i, i\in\mathcal{I}$, we have
\begin{enumerate}
\item The function $f^i_{j,\mu}(\mathbf{x})$ is convex in $\mathbf{x}^i$ and totally differentiable in $\mathbf{x}$.
\item The total gradient $\nabla f^i_{j,\mu}(\mathbf{x})$ is $\mathcal{L}$-Lipschitz continuous in $\mathbf{x}$, \textit{i.e.}, $\forall\mathbf{x},\mathbf{x}'\in\mathbb{R}^n$, $\|\nabla f^i_{j,\mu}(\mathbf{x}) - \nabla f^i_{j,\mu}(\mathbf{x}')\|\leq \mathcal{L}\|\mathbf{x} - \mathbf{x}'\|$; and satisfies $\|\nabla f^i_{j,\mu}(\mathbf{x}) - \nabla f^i_j(\mathbf{x})\|\leq\frac12(n+3)^{\frac32}\mathcal{L}\mu$.
\item The oracle $\pi^i_{jk}(\mathbf{x}_t)$ satisfies that $\mathbb{E}[\pi^i_{jk}(\mathbf{x}_t)|\mathcal{F}_t] = \nabla_{x^i_k} f^i_{j,\mu}(\mathbf{x}_t)$, and $\mathbb{E}[\|\pi^i_{jk}(\mathbf{x}_t)\|^2|\mathcal{F}_t] \leq 4(n+4)\|\nabla f^i_{j,\mu}(\mathbf{x}_t)\|^2+3(n+4)^3\mathcal{L}^2\mu^2$.
\end{enumerate}
\end{Lemma}

\subsection{Gaussian-Smoothed Game Formulation}

We define a Gaussian-smoothed $N$-coalition game, denoted by $\Gamma_\mu(\mathcal{I},\{f^i_\mu\},\{\mathbb{R}^{n_i}\})$, which has the same set of coalitions and action sets, but the cost function $f^i_{\mu}$ is given by
\begin{align*}
f^i_{\mu}(\mathbf{x}^i,\mathbf{x}^{-i}) \triangleq \frac1{n_i}\sum_{j=1}^{n_i}f^i_{j,\mu}(\mathbf{x}^i,\mathbf{x}^{-i}), \quad \forall i\in\mathcal{I},
\end{align*}
where $f^i_{j,\mu}$ is given by \eqref{eq:smooth_function_definition}.
The game mapping of $\Gamma_\mu$ is defined as $F_\mu(\mathbf{x}) \triangleq [\nabla_{\mathbf{x}^1}f^1_\mu(\mathbf{x})^\top,\ldots,\nabla_{\mathbf{x}^N}f^N_\mu(\mathbf{x})^\top]^\top$. 

In the next lemma, we proceed to show that $F_\mu(\mathbf{x})$ is strongly monotone, and characterize the distance between the NE of the smoothed game $\Gamma_\mu$ and the NE of the original game $\Gamma$ in terms of the smoothing parameter $\mu$.

\begin{Lemma}\label{lemma:NE_gap} Under Assumptions~\ref{assumption_local_f_lipschitz} and \ref{assumption_game_mapping}, for $\forall\mu\geq 0$, the smoothed game $\Gamma_\mu(\mathcal{I},\{f^i_\mu\},\{\mathbb{R}^{n_i}\})$ holds that
\begin{enumerate}
\item The game mapping $F_\mu(\mathbf{x})$ is $\chi$-strongly monotone.
\item It admits a unique NE (denoted by $\mathbf{x}_\mu^*$) satisfying
\begin{align*}
\|\mathbf{x}_\mu^* - \mathbf{x}^*\| \leq \frac{n(n+3)^{\frac32}\mathcal{L}\gamma}{2(1-\sqrt{1-\gamma\chi})}\mu,
\end{align*}
where $\mathbf{x}^*$ is the unique NE of the original game $\Gamma$, and $\gamma\in(0,\frac{\chi}{n^2\mathcal{L}^2}]$ is a constant.
\end{enumerate}
\end{Lemma}
\begin{Proof}
(1) Under Assumption~\ref{assumption_local_f_lipschitz}, we have $f^i_j(\mathbf{x})$ is differentiable in $\mathbf{x}^i$, and both $f^i_j(\mathbf{x})$ and $\nabla_{\mathbf{x}^i}f^i_j(\mathbf{x})$ are continuous in $\mathbf{x}$. Then, the integrand in \eqref{eq:smooth_function_definition} has both itself and its partial derivative with respect to $\mathbf{x}^i$ being continuous in $(\mathbf{x},\xi)$, Thus, for $\mu\geq 0$, we may derive $\nabla_{\mathbf{x}^i}f^i_{j,\mu}(\mathbf{x})$ by applying Leibniz integral rule to \eqref{eq:smooth_function_definition}, given by
\begin{align}
  \nabla_{\mathbf{x}^i}f^i_{j,\mu}(\mathbf{x}) = \frac1{\kappa}\int_{\mathbb{R}^n} \nabla_{\mathbf{x}^i}f^i_j(\mathbf{x}+\mu\xi)e^{-\frac12\|\xi\|^2}d\xi. \label{eq:smooth_gradient_definition}
\end{align}
For $\forall\mathbf{x},\mathbf{y}\in\mathbb{R}^n$, it can be shown that
\begin{align*}
&\langle F_\mu(\mathbf{x}) - F_\mu(\mathbf{y}), \mathbf{x} - \mathbf{y}\rangle = \sum_{i=1}^N\langle\nabla_{\mathbf{x}^i}f^i_\mu(\mathbf{x})-\nabla_{\mathbf{y}^i}f^i_\mu(\mathbf{y}),\\
&\mathbf{x}^i - \mathbf{y}^i\rangle=\frac1{\kappa}\int_{\mathbb{R}^n}\sum_{i=1}^N\langle\nabla_{\mathbf{x}^i}f^i(\mathbf{x}+\mu\xi)-\nabla_{\mathbf{y}^i}f^i(\mathbf{y}+\mu\xi),\\
&\mathbf{x}^i - \mathbf{y}^i\rangle e^{-\frac12\|\xi\|^2}d\xi=\frac1{\kappa}\int_{\mathbb{R}^n}\langle F(\mathbf{x}+\mu\xi)-F(\mathbf{y}+\mu\xi),\\
&(\mathbf{x}+\mu\xi) - (\mathbf{y}+\mu\xi)\rangle e^{-\frac12\|\xi\|^2}d\xi \geq \chi\|\mathbf{x} - \mathbf{y}\|^2,
\end{align*}
where the second equality is due to \eqref{eq:smooth_gradient_definition} and the last inequality follows from Assumption~\ref{assumption_game_mapping}.
Thus, we have that for $\mu\geq 0$, the game mapping of $\Gamma_\mu$ is $\chi$-strongly monotone. 

(2) With convexity of $f^i_\mu(\mathbf{x})$ in $\mathbf{x}^i$ (c.f. Lemma~\ref{lemma:property_f_mu}-(1)) and the strong monotonicity of $F_\mu(\mathbf{x})$, game $\Gamma_\mu$ admits a unique NE, denoted by $\mathbf{x}_\mu^*$. It is noted that an $N$-coalition game is equivalent to an $N$-player non-cooperative game with the cost being the coalition cost function, hence we have $F_\mu(\mathbf{x}^*_\mu) = \mathbf{0}_n$.

We use the notation $\bm{F}(\mathbf{x},\mu):\mathbb{R}^n\times\mathbb{R}_{\geq0}\to\mathbb{R}^n$ for the game mapping $F_\mu(\mathbf{x})$ of the smoothed game, and $\mathbf{x}^*(\mu)$ for the unique NE $\mathbf{x}_\mu^*$ to explicitly quantify the effect of the smoothing parameter $\mu$. For $\forall \mu\geq0$, it follows from part (1) that $\bm{F}(\mathbf{x},\mu)$ is $\chi$-strongly monotone in $\mathbf{x}$, \textit{i.e.}, $\forall \mathbf{x},\mathbf{x}'\in\mathbb{R}^n$,
\begin{align} \label{eq:map_strong_monotone}
\langle \mathbf{x}-\mathbf{x}', \bm{F}(\mathbf{x},\mu) - \bm{F}(\mathbf{x}',\mu)\rangle \geq \chi\|\mathbf{x}-\mathbf{x}'\|^2.
\end{align}
When $\mu = 0$, $f^i_{j,0}(\mathbf{x})$ reduces to $f^i_j(\mathbf{x})$ by definition~\eqref{eq:smooth_function_definition}. That implies the two games $\Gamma_\mu$ and $\Gamma$ are equal at $\mu=0$. Hence their respective unique NE are equal at $\mu=0$, \textit{i.e.}, $\mathbf{x}^*(0)= \mathbf{x}^*$. Now, establishing an upper bound for $\|\mathbf{x}_\mu^* - \mathbf{x}^*\|$ is equivalent to studying $\|\mathbf{x}^*(\mu) - \mathbf{x}^*(0)\|$, which is the Lipschitz property of $\mathbf{x}^*(\mu)$ at $\mu=0$.

The rest of the proof follows the idea in \cite[Thm.~2.1]{Dafermos1988}. We first show that $\bm{F}(\mathbf{x},\mu)$ is (i) Lipschitz continuous in $\mathbf{x}$ for $\forall \mu\geq0$, and (ii) Lipschitz continuous in $\mu$ at $\mu=0$ for $\forall \mathbf{x}\in\mathbb{R}^n$. For (i), it follows from Lemma~\ref{lemma:property_f_mu}-(2) that, $\forall \mathbf{x},\mathbf{x}'\in\mathbb{R}^n$, $\|\nabla_{\mathbf{x}^i}f^i_\mu(\mathbf{x})-\nabla_{\mathbf{x}^i}f^i_\mu(\mathbf{x}')\|=\|\sum_{j=1}^{n_i}(\nabla_{\mathbf{x}^i}f^i_{j,\mu}(\mathbf{x})-\nabla_{\mathbf{x}^i}f^i_{j,\mu}(\mathbf{x}'))\|\leq\sum_{j=1}^{n_i}\|\nabla f^i_{j,\mu}(\mathbf{x})-\nabla f^i_{j,\mu}(\mathbf{x}')\|\leq n_i\mathcal{L}\|\mathbf{x}-\mathbf{x}'\|$. Hence, we have $\|\bm{F}(\mathbf{x},\mu) - \bm{F}(\mathbf{x}',\mu)\| = (\sum_{i=1}^N\|\nabla_{\mathbf{x}^i}f^i_\mu(\mathbf{x})-\nabla_{\mathbf{x}^i}f^i_\mu(\mathbf{x}')\|^2)^{\frac12}$, which gives
\begin{align} \label{eq:map_lipschitz_x}
\|\bm{F}(\mathbf{x},\mu) - \bm{F}(\mathbf{x}',\mu)\| \leq n\mathcal{L}\|\mathbf{x}-\mathbf{x}'\|, \forall \mu\geq0.
\end{align}
For (ii), it follows from Lemma~\ref{lemma:property_f_mu}-(2) that, $\forall \mathbf{x}\in\mathbb{R}^n$, $\|\nabla_{\mathbf{x}^i}f^i_\mu(\mathbf{x})-\nabla_{\mathbf{x}^i}f^i(\mathbf{x})\|=\|\sum_{j=1}^{n_i}(\nabla_{\mathbf{x}^i}f^i_{j,\mu}(\mathbf{x})-\nabla_{\mathbf{x}^i}f^i_j(\mathbf{x}))\|\leq\sum_{j=1}^{n_i}\|\nabla f^i_{j,\mu}(\mathbf{x})-\nabla f^i_j(\mathbf{x})\|\leq \frac12n_i(n+3)^{\frac32}\mathcal{L}\mu$. Hence, we have
\begin{align} \label{eq:map_lipschitz_mu}
\|\bm{F}(\mathbf{x},\mu) - \bm{F}(\mathbf{x},0)\| \leq \frac12n(n+3)^{\frac32}\mathcal{L}\mu, \forall \mathbf{x}\in\mathbb{R}^n.
\end{align}
Now, we consider the map $M(\mathbf{x},\mu)\triangleq \mathbf{x}-\gamma\bm{F}(\mathbf{x},\mu)$. For $0<\gamma\leq\frac{\chi}{n^2\mathcal{L}^2}$, $\forall \mathbf{x},\mathbf{y}\in\mathbb{R}^n$,
\begin{align} 
&\|M(\mathbf{x},\mu) - M(\mathbf{y},\mu)\|^2 \leq \|\mathbf{x}-\mathbf{y}\|^2-2\gamma\langle \mathbf{x}-\mathbf{y}, \bm{F}(\mathbf{x},\mu) \nonumber\\
&\quad- \bm{F}(\mathbf{y},\mu) \rangle +\gamma^2\|\bm{F}(\mathbf{x},\mu) - \bm{F}(\mathbf{y},\mu)\|^2 \leq(1-2\gamma\chi\nonumber\\
&\quad+n^2\gamma^2\mathcal{L}^2)\|\mathbf{x}-\mathbf{y}\|^2\leq(1-\gamma\chi)\|\mathbf{x}-\mathbf{y}\|^2, \label{eq:map_contraction}
\end{align}
where the second inequality follows from \eqref{eq:map_strong_monotone} and \eqref{eq:map_lipschitz_x}. Hence, \eqref{eq:map_contraction} implies that the map $M(\mathbf{x},\mu)$ is a contraction with respect to $\mathbf{x}$. By the Banach fixed point theorem, the map $M(\mathbf{x},\mu)$ has a unique fixed point. On the other hand, any fixed point of the map $M(\mathbf{x},\mu)$ is a NE of game $\Gamma_\mu$ (since $\bm{F}(\mathbf{x},\mu) = F_\mu(\mathbf{x})=\mathbf{0}_n$). Thus, we deduce that the unique NE $\mathbf{x}^*(\mu)$ is the unique fixed point of the map $M(\mathbf{x},\mu)$. Then, it follows from \eqref{eq:map_contraction} that
\begin{align*}
\|\mathbf{x}^*(\mu) - \mathbf{x}^*(0)\| &= \|M(\mathbf{x}^*(\mu),\mu) - M(\mathbf{x}^*(0),0)\|\\
&\leq\|M(\mathbf{x}^*(\mu),\mu) - M(\mathbf{x}^*(0),\mu)\|\\
&\quad+\|M(\mathbf{x}^*(0),\mu) - M(\mathbf{x}^*(0),0)\|\\
&\leq\sqrt{1-\gamma\chi}\|\mathbf{x}^*(\mu) - \mathbf{x}^*(0)\|\\
&\quad+\|M(\mathbf{x}^*(0),\mu) - M(\mathbf{x}^*(0),0)\|.
\end{align*}
Noting that the last term
\begin{align*}
&\|M(\mathbf{x}^*(0),\mu) - M(\mathbf{x}^*(0),0)\| \\
= &\|(\mathbf{x}^*(0) - \gamma\bm{F}(\mathbf{x}^*(0),\mu)) - (\mathbf{x}^*(0) - \gamma\bm{F}(\mathbf{x}^*(0),0))\|\\
=&\gamma\|\bm{F}(\mathbf{x}^*(0),\mu)-\bm{F}(\mathbf{x}^*(0),0)\|\leq\frac12n(n+3)^{\frac32}\mathcal{L}\gamma\mu,
\end{align*}
where the last inequality is due to \eqref{eq:map_lipschitz_mu}. Combining the above two relations gives the desired result.
\end{Proof}

\section{Nash Equilibrium Seeking Strategy in $N$-Coalition Games}\label{sec:distr_opt}

In this section, we present the details of the NE seeking strategy. 

At time-step $t$, the player $j\in\mathcal{V}^i$ in each coalition $i\in\mathcal{I}$ needs to maintain the following variables: the player's own action variable $x^i_{j,t}$, auxiliary action variables $y^i_{jk,t}$, and gradient tracker variables $\phi^i_{jk,t}$ for $\forall k\in\mathcal{V}^i$. A gradient estimator $\pi^i_{jk}$ for $\forall k\in\mathcal{V}^i$ is also constructed to estimate the partial gradient based on the values of the local cost functions. The algorithm is initialized with arbitrary $x^i_{j,0},y^i_{jk,0}\in\mathbb{R}$ and $\phi^i_{jk,0} = \pi^i_{jk}(\mathbf{x}_0)$. Then, each player $j\in\mathcal{V}^i, i\in\mathcal{I}$ updates these variables according to the following update laws
\begin{subequations}\label{eq:algorithm}
\begin{align}
y^i_{jk,t+1} &= \sum_{l=1}^{n_i}[A^i]_{jl}y^i_{lk,t}-\alpha \phi^i_{jk,t},\label{eq:update_y}\\
x^i_{j,t+1} &= y^i_{jj,t+1}, \label{eq:update_x}\\
\phi^i_{jk,t+1}& = \sum_{l=1}^{n_i}[A^i]_{jl} \phi^i_{lk,t} + \pi^i_{jk}(\mathbf{x}_{t+1}) - \pi^i_{jk}(\mathbf{x}_t),\label{eq:update_phi}
\end{align}
\end{subequations}
where $\pi^i_{jk}(\mathbf{x}_t)$ is the gradient estimator given by \eqref{grad_oracle},
and $\alpha > 0$ is a constant step-size sequence. Thus, at time-step $t$, agent $j\in\mathcal{V}^i$ of coalition $i\in\mathcal{I}$ first update the auxiliary action variables $y^i_{jk,t}$, its decision variable $x^i_{j,t}$ and gradient tracker variables $\phi^i_{jk,t}, k\in\mathcal{V}^i$ according to the update laws \eqref{eq:algorithm}, and pass the updated auxiliary action variables $y^i_{jk,t}$ and gradient tracker variables $\phi^i_{jk,t}$, $k\in\mathcal{V}^i$ to its out-neighbors in the same coalition. Then the iteration continues. These procedures are tabulated in Algorithm~\ref{algo:rgf_d_dgd}.

\begin{algorithm}
\caption{NE seeking in $N$-coalition games}\label{algo:rgf_d_dgd}
\begin{algorithmic}[1]
\Initialize {set $x^i_{j,0},y^i_{jk,0}\in\mathbb{R}$\\ generate $\{\xi^i_{j,t}\}_{t\geq0}\sim\mathcal{N}(\mathbf{0}_n,I_n)$\\ set $\phi^i_{jk,0} = \pi^i_{jk}(\mathbf{x}_0)$, $\forall k\in\mathcal{V}^i$}
\Iteration {(communicate with neighbor $l$ within coalition $i$ to obtain $y^i_{lk,t-1}$, $\forall k\in\mathcal{V}^i$)\\update $y^i_{jk,t}$ based on \eqref{eq:update_y}\\ update $x^i_{j,t}$ based on \eqref{eq:update_x}\\($\mathbf{x}_t$ is broadcast to all players across all coalitions)\\compute $\pi^i_{jk}(\mathbf{x}_{t})$ based on \eqref{grad_oracle}, $\forall k\in\mathcal{V}^i$\\(communicate with neighbor $l$ within coalition $i$ to obtain $\phi^i_{lk,t-1}$, $\forall k\in\mathcal{V}^i$)\\update $\phi^i_{jk,t}$ based on \eqref{eq:update_phi}, $\forall k\in\mathcal{V}^i$}
\Output {$x^i_{j,t} \to x^{i*}_j$}
\end{algorithmic}
\end{algorithm}


\section{Convergence Analysis} \label{sec:conv_analysis}

To facilitate the convergence analysis, we first make some notations for easy presentation. Denote that $n_s\triangleq \sum_{i=1}^{N}n_i^2$ and $n_c\triangleq \sum_{i=1}^{N}n_i^3$. For $\forall k\in\mathcal{V}^i, i \in\mathcal{I}$, define that
\textcolor{black}{$\mathbf{y}^i_{k,t}\triangleq [y^i_{1k,t}, \ldots, y^i_{n_ik,t}]^\top\in\mathbb{R}^{n_i}$,
$\bar{y}^i_{k,t}\triangleq \frac1{n_i}\mathbf{1}_{n_i}^\top\mathbf{y}^i_{k,t}\in\mathbb{R}$,
$\bar{\mathbf{y}}^i_t\triangleq [\bar{y}^i_{1,t},\ldots,\bar{y}^i_{n_i,t}]^\top\in\mathbb{R}^{n_i}$,
$\bar{\mathbf{y}}_t\triangleq [\bar{\mathbf{y}}^{1\top}_t,\ldots,\bar{\mathbf{y}}^{N\top}_t]^\top\in\mathbb{R}^n$,}
$\phi_t \triangleq [\phi^1_{11,t},\ldots,\phi^N_{n_Nn_N,t}]^\top\in\mathbb{R}^n$,
$\phi^i_{k,t} \triangleq [\phi^i_{1k,t},\ldots,\phi^i_{n_ik,t}]^\top\in\mathbb{R}^{n_i}$, 
$\bar{\phi}^i_{k,t} \triangleq \frac1{n_i}\mathbf{1}_{n_i}^\top\phi^i_{k,t}\in\mathbb{R}$,
\textcolor{black}{$\bar{\bm{\phi}}^i_t \triangleq [\bar{\phi}^i_{1,t},\ldots,\bar{\phi}^i_{n_i,t}]^\top\in\mathbb{R}^{n_i}$,
$\bar{\bm{\phi}}_t \triangleq [\bar{\bm{\phi}}^{1\top}_t,\ldots,\bar{\bm{\phi}}^{N\top}_t]^\top\in\mathbb{R}^n$,}
$\pi^i_{k} \triangleq [\pi^i_{1k},\ldots,\pi^i_{n_ik}]^\top\in\mathbb{R}^{n_i}$, $\bar{\pi}^i_k \triangleq \frac1{n_i}\mathbf{1}_{n_i}^\top\pi^i_k\in\mathbb{R}$, and
$\nabla_{x^i_k} \mathbf{f}^i_\mu(\mathbf{x}) \triangleq [\nabla_{x^i_k} f^i_{1,\mu}(\mathbf{x}),\ldots,\nabla_{x^i_k} f^i_{n_i,\mu}(\mathbf{x})]^\top\in\mathbb{R}^{n_i}$.
Then, the update laws \eqref{eq:update_y} and \eqref{eq:update_phi} can be compactly written as
\begin{subequations}
\begin{align}
\color{black}\mathbf{y}^i_{k,t+1} &\color{black}= A^i\mathbf{y}^i_{k,t}-\alpha \phi^i_{k,t},\label{eq:update_y_compact}\\
\phi^i_{k,t+1} &= A^i \phi^i_{k,t} + \pi^i_k(\mathbf{x}_{t+1}) - \pi^i_k(\mathbf{x}_t).\label{eq:update_phi_compact}
\end{align}
\end{subequations}
\textcolor{black}{Multiplying $\frac1{n_i}\mathbf{1}_{n_i}^\top$ from the left on both sides of \eqref{eq:update_y_compact} and augmenting the relation into a compact form, we obtain}
\begin{align}
\color{black}\bar{\mathbf{y}}_{t+1} = \bar{\mathbf{y}}_t - \alpha \bar{\bm{\phi}}_t, \label{eq:update_y_bar_compact}
\end{align}

Now, we are ready for the convergence analysis of the proposed algorithm, which consists of two parts: 1) in Sec.~\ref{subsec:auxiliary_results}, we derive the inequality iterations of three major quantities: i) \textcolor{black}{$\sum_{i=1}^N\sum_{k=1}^{n_i}\mathbb{E}[\|\mathbf{y}^i_{k,t}-\mathbf{1}_{n_i}\bar{y}^i_{k,t}\|^2]$, the total consensus error of the auxiliary variables, ii) $\mathbb{E}[\|\bar{\mathbf{y}}_t-\mathbf{x}^*_\mu\|^2]$, the gap between the stacked averaged auxiliary variable and the NE of game $\Gamma_\mu$,} and iii) $\sum_{i=1}^N\sum_{k=1}^{n_i}\mathbb{E}[\|\phi^i_{k,t} - \mathbf{1}_{n_i}\bar{\phi}^i_{k,t}\|^2]$, the total gradient tracking error; 2) in Sec.~\ref{subsec:main_results}, we establish a linear system of these three inequalities to analyze its convergence, \textcolor{black}{and characterize the gap between the players' actions and the NE of game $\Gamma_\mu$, $\mathbb{E}[\|\mathbf{x}_t-\mathbf{x}^*_\mu\|^2]$ in terms of the aforementioned quantities, and finally deduce the gap between the players' actions and the NE of game $\Gamma$, $\mathbb{E}[\|\mathbf{x}_t-\mathbf{x}^*\|^2]$ with the obtained results in Lemma~\ref{lemma:NE_gap}-(2).}

\subsection{Auxiliary Results}\label{subsec:auxiliary_results}
First, we derive some basic properties on the averaged gradient tracker $\bar{\phi}^i_{k,t}$, and provide a bound on the stacked gradient tracker $\phi^i_{k,t}$ in Lemmas~\ref{lemma:averaged_grad_tracker} and \ref{lemma:stacked_grad_tracker}, respectively.

\begin{Lemma}\label{lemma:averaged_grad_tracker}
Under Assumptions~\ref{assumption_graph}, \ref{assumption_local_f_lipschitz} and \ref{assumption_game_mapping}, the averaged gradient tracker $\bar{\phi}^i_{k,t}, \forall k\in\mathcal{V}^i, i \in\mathcal{I}$ holds that
\begin{enumerate}
\item $\bar{\phi}^i_{k,t}=\bar{\pi}^i_{k}(\mathbf{x}_t)$,
\item $\mathbb{E}[\bar{\phi}^i_{k,t}|\mathcal{F}_t] = \nabla_{x^i_k} f^i_\mu(\mathbf{x}_t)$,
\item \textcolor{black}{$\mathbb{E}[\|\bar{\phi}^i_{k,t}\|^2|\mathcal{F}_t] \leq 12(n+4)\mathcal{L}^2\sum_{i=1}^N\sum_{k=1}^{n_i}\|\mathbf{y}^i_{k,t}-\mathbf{1}_{n_i}\bar{y}^i_{k,t}\|^2+12(n+4)\mathcal{L}^2\|\bar{\mathbf{y}}_t-\mathbf{x}^*_\mu\|^2+12(n+4)G^2+3(n+4)^3\mathcal{L}^2\mu^2$, where $G\triangleq \max_{j\in\mathcal{V}^i,i\in\mathcal{I}}\|\nabla f^i_{j,\mu}(\mathbf{x}^*_\mu)\|$.}
\end{enumerate}
\end{Lemma}
\begin{Proof}
For (1), multiplying $\frac1{n_i}\mathbf{1}_{n_i}^\top$ from the left on both sides of \eqref{eq:update_phi_compact}, and noting that $A^i$ is doubly stochastic, we have
$\bar{\phi}^i_{k,t+1} = \bar{\phi}^i_{k,t} + \bar{\pi}^i_k(\mathbf{x}_{t+1}) - \bar{\pi}^i_k(\mathbf{x}_t)$.
Recursively expanding the above relation and noting that $\phi^i_{k,0} = \pi^i_k(\mathbf{x}_0)$ completes the proof.

For (2), following the result of part (1) and Lemma~\ref{lemma:property_f_mu}-3), we obtain
$\mathbb{E}[\bar{\phi}^i_{k,t}|\mathcal{F}_t] = \mathbb{E}[\bar{\pi}^i_k(\mathbf{x}_t)|\mathcal{F}_t]=\frac1{n_i}\mathbf{1}_{n_i}^\top\mathbb{E}[\pi^i_k|\mathcal{F}_t]=\frac1{n_i}\mathbf{1}_{n_i}^\top\nabla_{x^i_k} \mathbf{f}^i_\mu(\mathbf{x}_t)=\nabla_{x^i_k} f^i_\mu(\mathbf{x}_t)$.

For (3), It is noted that
\begin{align*}
&\mathbb{E}[\|\bar{\phi}^i_{k,t}\|^2|\mathcal{F}_t]=\mathbb{E}[\|\bar{\pi}^i_{k}(\mathbf{x}_t)\|^2|\mathcal{F}_t] \\
&\quad= \frac1{n_i^2}\mathbb{E}[\|\mathbf{1}_{n_i}^\top\pi^i_{k}(\mathbf{x}_t)\|^2|\mathcal{F}_t]\leq \frac1{n_i}\sum_{j=1}^{n_i}\mathbb{E}[\|\pi^i_{jk}(\mathbf{x}_t)\|^2|\mathcal{F}_t].
\end{align*}
\color{black}
From Lemma~\ref{lemma:property_f_mu}-(3),
\begin{align}
&\mathbb{E}[\|\pi^i_{jk}(\mathbf{x}_t)\|^2|\mathcal{F}_t]\leq 4(n+4)\|\nabla f^i_{j,\mu}(\mathbf{x}_t)\|^2+3(n+4)^3\mathcal{L}^2\mu^2\nonumber\\
&\leq 12(n+4)\|\nabla f^i_{j,\mu}(\mathbf{x}_t)-\nabla f^i_{j,\mu}(\bar{\mathbf{y}}_t)\|^2\nonumber\\
&\quad+12(n+4)\|\nabla f^i_{j,\mu}(\bar{\mathbf{y}}_t)-\nabla f^i_{j,\mu}(\mathbf{x}^*_\mu)\|^2\nonumber\\
&\quad+12(n+4)\|\nabla f^i_{j,\mu}(\mathbf{x}^*_\mu)\|^2+3(n+4)^3\mathcal{L}^2\mu^2\nonumber\\
&\leq12(n+4)\mathcal{L}^2\|\mathbf{x}_t-\bar{\mathbf{y}}_t\|^2+12(n+4)\mathcal{L}^2\|\bar{\mathbf{y}}_t-\mathbf{x}^*_\mu\|^2\nonumber\\
&\quad+12(n+4)G^2+3(n+4)^3\mathcal{L}^2\mu^2\nonumber\\
&\leq12(n+4)\mathcal{L}^2\sum_{i=1}^N\sum_{k=1}^{n_i}\|\mathbf{y}^i_{k,t}-\mathbf{1}_{n_i}\bar{y}^i_{k,t}\|^2+3(n+4)^3\mathcal{L}^2\mu^2\nonumber\\
&\quad+12(n+4)\mathcal{L}^2\|\bar{\mathbf{y}}_t-\mathbf{x}^*_\mu\|^2+12(n+4)G^2,\label{eq:pi_i_jk_bound}
\end{align}
where $G\triangleq \max_{j\in\mathcal{V}^i,i\in\mathcal{I}}\|\nabla f^i_{j,\mu}(\mathbf{x}^*_\mu)\|$ and the last inequality follows from \eqref{eq:update_x} that
\begin{align}
\|\mathbf{x}_t-\bar{\mathbf{y}}_t\|^2&=\sum_{i=1}^N\sum_{k=1}^{n_i}\|y^i_{kk,t}-\bar{y}^i_{k,t}\|^2 \nonumber\\
&\leq\sum_{i=1}^N\sum_{k=1}^{n_i}\|\mathbf{y}^i_{k,t}-\mathbf{1}_{n_i}\bar{y}^i_{k,t}\|^2. \label{eq:x_minus_ybar}
\end{align}
The proof is completed by combining the preceding relations.
\end{Proof}

\begin{Lemma}\label{lemma:stacked_grad_tracker}
Under Assumptions~\ref{assumption_graph}, \ref{assumption_local_f_lipschitz} and \ref{assumption_game_mapping}, the stacked gradient tracker $\{\phi^i_{k,t}\}_{t\geq0}, \forall k\in\mathcal{V}^i, i \in\mathcal{I}$ holds that
\begin{align*}
&\mathbb{E}[\|\phi^i_{k,t}\|^2|\mathcal{F}_t]\leq 2\mathbb{E}[\|\phi^i_{k,t} - \mathbf{1}_{n_i}\bar{\phi}^i_{k,t}\|^2|\mathcal{F}_t]\\
&\quad\color{black}+ 24n_i^2(n+4)\mathcal{L}^2\sum_{i=1}^N\sum_{k=1}^{n_i}\|\mathbf{y}^i_{k,t}-\mathbf{1}_{n_i}\bar{y}^i_{k,t}\|^2\\
&\quad\color{black}+24n_i^2(n+4)\mathcal{L}^2\|\bar{\mathbf{y}}_t-\mathbf{x}^*_\mu\|^2\\
&\quad\color{black}+24n_i^2(n+4)G^2+6n_i^2(n+4)^3\mathcal{L}^2\mu^2.
\end{align*}
\end{Lemma}
\begin{Proof}
For $\forall k\in\mathcal{V}^i, i \in\mathcal{I}$, we have
\begin{align*}
\|\phi^i_{k,t}\|^2 &\leq 2\|\phi^i_{k,t} - \mathbf{1}_{n_i}\bar{\phi}^i_{k,t}\|^2 + 2n_i^2\|\bar{\phi}^i_{k,t}\|^2.
\end{align*}
The proof is completed by taking the conditional expectation on $\mathcal{F}_t$ and substituting Lemma~\ref{lemma:averaged_grad_tracker}-(3).
\end{Proof}

In the subsequent Lemmas~\ref{lemma:conensus}, \ref{lemma:optimality_gap} and \ref{lemma:grad_track_error}, we characterize the inequality iterations of the aforementioned three major quantities. Specifically, we begin with \textcolor{black}{the total consensus error of the auxiliary variables $\sum_{i=1}^N\sum_{k=1}^{n_i}\|\mathbf{y}^i_{k,t}-\mathbf{1}_{n_i}\bar{y}^i_{k,t}\|^2$,} in terms of the sum of the previous iterations of the three quantities, and some constants.

\color{black}
\begin{Lemma}\label{lemma:conensus}
Under Assumptions~\ref{assumption_graph}, \ref{assumption_local_f_lipschitz} and \ref{assumption_game_mapping}, the total consensus error of the auxiliary variables $\sum_{i=1}^N\sum_{k=1}^{n_i}\|\mathbf{y}^i_{k,t}-\mathbf{1}_{n_i}\bar{y}^i_{k,t}\|^2$ satisfies that
\begin{align*}
&\sum_{i=1}^N\sum_{k=1}^{n_i}\mathbb{E}[\|\mathbf{y}^i_{k,t+1}-\mathbf{1}_{n_i}\bar{y}^i_{k,t+1}\|^2|\mathcal{F}_t]\leq\bigg(\frac{1+\bar{\sigma}^2}2 \\
&\quad+ 24(n+4)n_c\varsigma\mathcal{L}^2\alpha^2\bigg)\sum_{i=1}^N\sum_{k=1}^{n_i}\|\mathbf{y}^i_{k,t}-\mathbf{1}_{n_i}\bar{y}^i_{k,t}\|^2\\
&\quad+24(n+4)n_c\varsigma\mathcal{L}^2\alpha^2\|\bar{\mathbf{y}}_t-\mathbf{x}^*_\mu\|^2\\
&\quad+2\varsigma\alpha^2\sum_{i=1}^N\sum_{k=1}^{n_i}\mathbb{E}[\|\phi^i_{k,t} - \mathbf{1}_{n_i}\bar{\phi}^i_{k,t}\|^2|\mathcal{F}_t]\\
&\quad+24(n+4)n_c\varsigma G^2\alpha^2+6(n+4)^3n_c\varsigma\mathcal{L}^2\mu^2\alpha^2.
\end{align*}
\end{Lemma}
\begin{Proof}
It follows from \eqref{eq:update_y_compact} that for $i\in\mathcal{I}$
\begin{align*}
&\|\mathbf{y}^i_{k,t+1}-\mathbf{1}_{n_i}\bar{y}^i_{k,t+1}\|^2 = \bigg\|A^i\mathbf{y}^i_{k,t} -\alpha\phi^i_{k,t}\\
&\quad\quad-\frac1{n_i}\mathbf{1}_{n_i}\mathbf{1}_{n_i}^\top (A^i\mathbf{y}^i_{k,t} -\alpha\phi^i_{k,t})\bigg\|^2\\
&\quad\leq\|A^i\mathbf{y}^i_{k,t}-\mathbf{1}_{n_i}\bar{y}^i_{k,t}\|^2 + \alpha^2\bigg\|\bigg(I_{n_i}-\frac1{n_i}\mathbf{1}_{n_i}\mathbf{1}_{n_i}^\top\bigg)\phi^i_{k,t}\bigg\|^2 \\
&\quad\quad-2\bigg\langle A^i\mathbf{y}^i_{k,t}-\mathbf{1}_{n_i}\bar{y}^i_{k,t}, \alpha\bigg(I_{n_i}-\frac1{n_i}\mathbf{1}_{n_i}\mathbf{1}_{n_i}^\top\bigg)\phi^i_{k,t}\bigg\rangle,
\end{align*}
Taking the conditional expectation on $\mathcal{F}_t$ and noting that $\|I_{n_i}-\frac1{n_i}\mathbf{1}_{n_i}\mathbf{1}_{n_i}^\top\| = 1$, we obtain
\begin{align*}
&\mathbb{E}[\|\mathbf{y}^i_{k,t+1}-\mathbf{1}_{n_i}\bar{y}^i_{k,t+1}\|^2|\mathcal{F}_t] \leq\sigma_{A^i}^2\|\mathbf{y}^i_{k,t}-\mathbf{1}_{n_i}\bar{y}^i_{k,t}\|^2\\
&+ \alpha^2\mathbb{E}[\|\phi^i_{k,t}\|^2|\mathcal{F}_t] +\frac{1-\sigma_{A^i}^2}{2\sigma_{A^i}^2}\mathbb{E}[\|A^i\mathbf{y}^i_{k,t}-\mathbf{1}_{n_i}\bar{y}^i_{k,t}\|^2|\mathcal{F}_t]\\
&+\frac{2\sigma_{A^i}^2}{1-\sigma_{A^i}^2}\alpha^2\mathbb{E}[\|\phi^i_{k,t}\|^2|\mathcal{F}_t]\leq\frac{1+\sigma_{A^i}^2}{1-\sigma_{A^i}^2}\alpha^2\mathbb{E}[\|\phi^i_{k,t}\|^2|\mathcal{F}_t]\\
&+\frac{1+\sigma_{A^i}^2}2\mathbb{E}[\|\mathbf{y}^i_{k,t}-\mathbf{1}_{n_i}\bar{y}^i_{k,t}\|^2|\mathcal{F}_t]\leq\varsigma\alpha^2\mathbb{E}[\|\phi^i_{k,t}\|^2|\mathcal{F}_t]\\
&+\frac{1+\bar{\sigma}^2}{2}\|\mathbf{y}^i_{k,t}-\mathbf{1}_{n_i}\bar{y}^i_{k,t}\|^2.
\end{align*}
Applying Lemma~\ref{lemma:stacked_grad_tracker} and summing over $k=1$ to $n_i$, $i=1$ to $N$ complete the proof.
\end{Proof}
\color{black}

Then, we quantify the gap \textcolor{black}{$\|\bar{\mathbf{y}}_t-\mathbf{x}^*_\mu\|^2$}, in terms of the combination of the previous iterations \textcolor{black}{of itself and the total consensus error of the auxiliary variables}, and some constants.

\begin{Lemma}\label{lemma:optimality_gap}
Under Assumptions~\ref{assumption_graph}, \ref{assumption_local_f_lipschitz} and \ref{assumption_game_mapping}, the gap between \textcolor{black}{the stacked averaged auxiliary variable and the NE of game $\Gamma_\mu$, $\|\bar{\mathbf{y}}_t-\mathbf{x}^*_\mu\|^2$} satisfies that
\begin{align*}
&\color{black}\mathbb{E}[\|\bar{\mathbf{y}}_{t+1}-\mathbf{x}^*_\mu\|^2|\mathcal{F}_t]\leq(1-\chi\alpha+12n(n+4)\mathcal{L}^2\alpha^2)\|\bar{\mathbf{y}}_t \\
&\quad\color{black}- \mathbf{x}^*_\mu\|^2+\bigg(\frac{\alpha\mathcal{L}^2}{\chi}+12n(n+4)\mathcal{L}^2\alpha^2\bigg)\sum_{i=1}^N\sum_{k=1}^{n_i}\|\mathbf{y}^i_{k,t}\\
&\quad\color{black}-\mathbf{1}_{n_i}\bar{y}^i_{k,t}\|^2+3n(n+4)^3\mathcal{L}^2\mu^2\alpha^2+12n(n+4)G^2\alpha^2.
\end{align*}
\end{Lemma}
\begin{Proof}
It follows from \textcolor{black}{\eqref{eq:update_y_bar_compact}} that
\begin{align*}
\color{black}\bar{\mathbf{y}}_{t+1}-\mathbf{x}^*_\mu = \bar{\mathbf{y}}_t -\mathbf{x}^*_\mu - \alpha \bar{\bm{\phi}}_t.
\end{align*}
Taking the Euclidean norm on both sides and noting that \textcolor{black}{$F_\mu(\mathbf{x}^*_\mu)=\mathbf{0}_n$,}
\begin{subequations}\label{eq:optimality_gap_expand}
\begin{align}
&\color{black}\|\bar{\mathbf{y}}_{t+1}-\mathbf{x}^*_\mu\|^2 \leq \|\bar{\mathbf{y}}_t -\mathbf{x}^*_\mu\|^2 + \alpha^2\sum_{i=1}^N\sum_{k=1}^{n_i} \|\bar{\phi}^i_{k,t}\|^2\nonumber\\
&\quad\color{black}-2\alpha\sum_{i=1}^N\sum_{k=1}^{n_i}\langle \bar{y}^i_{k,t} - x^{i*}_{k,\mu}, \bar{\phi}^i_{k,t} -\nabla_{x^i_k} f^i_\mu(\mathbf{x}_t) \rangle\label{eq:optimality_gap_expand_1}\\
&\quad\color{black}-2\alpha\langle \bar{\mathbf{y}}_t -\mathbf{x}^*_\mu, F_\mu(\mathbf{x}_t) -  F_\mu(\bar{\mathbf{y}}_t)\rangle\label{eq:optimality_gap_expand_2}\\
&\quad\color{black}-2\alpha\langle \bar{\mathbf{y}}_t -\mathbf{x}^*_\mu, F_\mu(\bar{\mathbf{y}}_t) - F_\mu(\mathbf{x}^*_\mu)\rangle. \label{eq:optimality_gap_expand_3}
\end{align}
\end{subequations}
For \eqref{eq:optimality_gap_expand_1}, it follows from Lemma~\ref{lemma:averaged_grad_tracker}-(2) that
\begin{align*}
&\mathbb{E}\bigg[-2\alpha\sum_{i=1}^N\sum_{k=1}^{n_i}\langle \bar{y}^i_{k,t} - x^{i*}_{k,\mu}, \bar{\phi}^i_{k,t} - \nabla_{x^i_k} f^i_\mu(\mathbf{x}_t)\rangle\bigg|\mathcal{F}_t\bigg]=0.
\end{align*}
For \eqref{eq:optimality_gap_expand_2}, it follows from Lemma~\ref{lemma:property_f_mu}-(2) that
\begin{align*}
&\color{black}-2\alpha\langle \bar{\mathbf{y}}_t -\mathbf{x}^*_\mu, F_\mu(\mathbf{x}_t) -  F_\mu(\bar{\mathbf{y}}_t)\rangle\\
&\color{black}\leq2\alpha\sum_{i=1}^N\sum_{k=1}^{n_i}\|\bar{y}^i_{k,t} - x^{i*}_{k,\mu}\|\|\nabla_{x^i_k} f^i_\mu(\mathbf{x}_t) - \nabla_{x^i_k} f^i_\mu(\bar{\mathbf{y}}_t)\|\\
&\color{black}\leq2\alpha\sum_{i=1}^N\sum_{k=1}^{n_i}\|\bar{y}^i_{k,t} - x^{i*}_{k,\mu}\|\mathcal{L}\|\mathbf{x}_t-\bar{\mathbf{y}}_t\|\\
&\color{black}\leq\chi\alpha\|\bar{\mathbf{y}}_t -\mathbf{x}^*_\mu\|^2+\frac{\alpha\mathcal{L}^2}\chi\|\mathbf{x}_t-\bar{\mathbf{y}}_t\|^2\\
&\color{black}\leq\chi\alpha\|\bar{\mathbf{y}}_t -\mathbf{x}^*_\mu\|^2+\frac{\alpha\mathcal{L}^2}\chi\sum_{i=1}^N\sum_{k=1}^{n_i}\|\mathbf{y}^i_{k,t}-\mathbf{1}_{n_i}\bar{y}^i_{k,t}\|^2,
\end{align*}
\textcolor{black}{where the last inequality follows from \eqref{eq:x_minus_ybar}.
For \eqref{eq:optimality_gap_expand_3}, it follows from Lemma~\ref{lemma:NE_gap}-(1) that}
\begin{align*}
\color{black}-2\alpha\langle \bar{\mathbf{y}}_t -\mathbf{x}^*_\mu, F_\mu(\bar{\mathbf{y}}_t) - F_\mu(\mathbf{x}^*_\mu)\rangle\leq-2\chi\alpha\|\bar{\mathbf{y}}_t -\mathbf{x}^*_\mu\|^2.
\end{align*}
Taking the conditional expectation on $\mathcal{F}_t$ for \eqref{eq:optimality_gap_expand}, and substituting the above three relations and Lemma~\ref{lemma:averaged_grad_tracker}-(3) yield the desired result.
\end{Proof}

Finally, we derive a bound on the total gradient tracking error $\sum_{i=1}^N\sum_{k=1}^{n_i}\|\phi^i_{k,t}- \mathbf{1}_{n_i}\bar{\phi}^i_{k,t}\|^2$, in terms of the combination of the previous iterations of the three quantities, and some constants. 

\begin{Lemma}\label{lemma:grad_track_error}
Under Assumptions~\ref{assumption_graph}, \ref{assumption_local_f_lipschitz} and \ref{assumption_game_mapping}, the total gradient tracking error $\sum_{i=1}^N\sum_{k=1}^{n_i}\|\phi^i_{k,t} - \mathbf{1}_{n_i}\bar{\phi}^i_{k,t}\|^2$ satisfies
\begin{align*}
&\color{black}\sum_{i=1}^N\sum_{k=1}^{n_i}\mathbb{E}[\|\phi^i_{k,t+1} - \mathbf{1}_{n_i}\bar{\phi}^i_{k,t+1}\|^2|\mathcal{F}_t]\leq\bigg(\frac{1+\bar{\sigma}^2}{2}\\
&\quad\color{black}+48(n+4)n_s\varsigma^2\mathcal{L}^2\alpha^2\bigg)\sum_{i=1}^N\sum_{k=1}^{n_i}\mathbb{E}[\|\phi^i_{k,t}- \mathbf{1}_{n_i}\bar{\phi}^i_{k,t}\|^2|\mathcal{F}_t] \\
&\color{black}+24(n+4)n_s\varsigma\mathcal{L}^2\bigg(\frac{3+\bar{\sigma}^2}2+\frac{\alpha\mathcal{L}^2}{\chi}+12n(n+4)\mathcal{L}^2\alpha^2 \\
&\quad\color{black}+ 24(n+4)n_c\varsigma\mathcal{L}^2\alpha^2\bigg)\sum_{i=1}^N\sum_{k=1}^{n_i}\|\mathbf{y}^i_{k,t}-\mathbf{1}_{n_i}\bar{y}^i_{k,t}\|^2\\
&\color{black}+24(n+4)n_s\varsigma\mathcal{L}^2[2+12n(n+4)\mathcal{L}^2\alpha^2\\
&\quad\color{black}+24(n+4)n_c\varsigma\mathcal{L}^2\alpha^2]\|\bar{\mathbf{y}}_t - \mathbf{x}^*_\mu\|^2\\
&\color{black}+24(n+4)n_s\varsigma\mathcal{L}^2[24(n+4)n_c\varsigma G^2\alpha^2\\
&\quad\color{black}+6(n+4)^3n_c\varsigma\mathcal{L}^2\mu^2\alpha^2+3n(n+4)^3\mathcal{L}^2\mu^2\alpha^2\\
&\quad\color{black}+12n(n+4)G^2\alpha^2]+48(n+4)n_s\varsigma G^2\\
&\quad\color{black}+12(n+4)^3n_s\varsigma\mathcal{L}^2\mu^2.
\end{align*}
\end{Lemma}
\begin{Proof}
It follows from \eqref{eq:update_phi_compact} that
\begin{align*}
&\|\phi^i_{k,t+1} - \mathbf{1}_{n_i}\bar{\phi}^i_{k,t+1}\|^2= \|A^i\phi^i_{k,t} - \mathbf{1}_{n_i}\bar{\phi}^i_{k,t}\|^2 \\
&\quad+ \bigg\|\bigg(I_{n_i}-\frac1{n_i}\mathbf{1}_{n_i}\mathbf{1}_{n_i}^\top\bigg)(\pi^i_k(\mathbf{x}_{t+1}) - \pi^i_k(\mathbf{x}_t))\bigg\|^2 \\
&\quad+2\bigg\langle A^i\phi^i_{k,t}- \mathbf{1}_{n_i}\bar{\phi}^i_{k,t},\\
&\quad\quad \bigg(I_{n_i}-\frac1{n_i}\mathbf{1}_{n_i}\mathbf{1}_{n_i}^\top\bigg)(\pi^i_k(\mathbf{x}_{t+1}) - \pi^i_k(\mathbf{x}_t) )\bigg\rangle.
\end{align*}
Taking the conditional expectation on $\mathcal{F}_t$, we have
\begin{align}
&\mathbb{E}[\|\phi^i_{k,t+1} - \mathbf{1}_{n_i}\bar{\phi}^i_{k,t+1}\|^2|\mathcal{F}_t] \leq \sigma_{A^i}^2\mathbb{E}[\|\phi^i_{k,t} \nonumber\\
&\quad\quad- \mathbf{1}_{n_i}\bar{\phi}^i_{k,t}\|^2|\mathcal{F}_t] + \mathbb{E}[\|\pi^i_k(\mathbf{x}_{t+1}) - \pi^i_k(\mathbf{x}_t)\|^2|\mathcal{F}_t]\nonumber\\
&\quad +2\mathbb{E}[\| A^i\phi^i_{k,t} - \mathbf{1}_{n_i}\bar{\phi}^i_{k,t}\| \|\pi^i_k(\mathbf{x}_{t+1}) - \pi^i_k(\mathbf{x}_t) \||\mathcal{F}_t]\nonumber\\
&\leq \sigma_{A^i}^2\mathbb{E}[\|\phi^i_{k,t} - \mathbf{1}_{n_i}\bar{\phi}^i_{k,t}\|^2|\mathcal{F}_t] + \mathbb{E}[\|\pi^i_k(\mathbf{x}_{t+1}) \nonumber\\
&\quad\quad- \pi^i_k(\mathbf{x}_t)\|^2|\mathcal{F}_t] +\frac{1-\sigma_{A^i}^2}2\mathbb{E}[\|\phi^i_{k,t} - \mathbf{1}_{n_i}\bar{\phi}^i_{k,t}\|^2|\mathcal{F}_t]\nonumber\\
&\quad+\frac{2\sigma_{A^i}^2}{1-\sigma_{A^i}^2}\mathbb{E}[\|\pi^i_k(\mathbf{x}_{t+1}) - \pi^i_k(\mathbf{x}_t)\|^2|\mathcal{F}_t]\nonumber\\
&\leq\frac{1+\bar{\sigma}^2}{2}\mathbb{E}[\|\phi^i_{k,t} - \mathbf{1}_{n_i}\bar{\phi}^i_{k,t}\|^2|\mathcal{F}_t]\nonumber\\
&\quad+\varsigma\mathbb{E}[\|\pi^i_k(\mathbf{x}_{t+1}) - \pi^i_k(\mathbf{x}_t)\|^2|\mathcal{F}_t],\label{eq:grad_track_error_term123}
\end{align}
\color{black}where the last term of \eqref{eq:grad_track_error_term123} follows from \eqref{eq:pi_i_jk_bound} that
\begin{align*}
&\mathbb{E}[\|\pi^i_k(\mathbf{x}_{t+1}) - \pi^i_k(\mathbf{x}_t)\|^2|\mathcal{F}_t]\\
&\quad\leq 2\sum_{j=1}^{n_i}\mathbb{E}[\|\pi^i_{jk}(\mathbf{x}_{t+1})\|^2|\mathcal{F}_t]+2\sum_{j=1}^{n_i}\mathbb{E}[\|\pi^i_{jk}(\mathbf{x}_t)\|^2|\mathcal{F}_t]\\
&\quad\leq 24n_i(n+4)\mathcal{L}^2\sum_{i=1}^N\sum_{k=1}^{n_i}\mathbb{E}[\|\mathbf{y}^i_{k,t+1}-\mathbf{1}_{n_i}\bar{y}^i_{k,t+1}\|^2|\mathcal{F}_t]\\
&\quad+24n_i(n+4)\mathcal{L}^2\mathbb{E}[\|\bar{\mathbf{y}}_{t+1}-\mathbf{x}^*_\mu\|^2|\mathcal{F}_t]\\
&\quad+24n_i(n+4)\mathcal{L}^2\sum_{i=1}^N\sum_{k=1}^{n_i}\|\mathbf{y}^i_{k,t}-\mathbf{1}_{n_i}\bar{y}^i_{k,t}\|^2\\
&\quad+24n_i(n+4)\mathcal{L}^2\|\bar{\mathbf{y}}_t-\mathbf{x}^*_\mu\|^2\\
&\quad+48n_i(n+4)G^2+12n_i(n+4)^3\mathcal{L}^2\mu^2.
\end{align*}
Applying Lemmas~\ref{lemma:conensus} and \ref{lemma:optimality_gap} in the above relation, and summing \eqref{eq:grad_track_error_term123} over $k=1$ to $n_i$, $i=1$ to $N$ complete the proof.\color{black}
\end{Proof}

\subsection{Main Results}\label{subsec:main_results}

Next, we are ready to analyze the convergence of the proposed algorithm with the following definitions.
\begin{align*}
\Psi_t &\triangleq \begin{bmatrix}\color{black}\begin{smallmatrix} \sum_{i=1}^N\sum_{k=1}^{n_i}\mathbb{E}[\|\mathbf{y}^i_{k,t}-\mathbf{1}_{n_i}\bar{y}^i_{k,t}\|^2]\\ \mathbb{E}[\|\bar{\mathbf{y}}_t - \mathbf{x}^*_\mu\|^2]\\\sum_{i=1}^N\sum_{k=1}^{n_i}\mathbb{E}[\|\phi^i_{k,t} - \mathbf{1}_{n_i}\bar{\phi}^i_{k,t}\|^2]\end{smallmatrix} \end{bmatrix}, \\
\mathbf{M}_\alpha &\triangleq\begin{bmatrix}\color{black}\begin{smallmatrix}1-k_1+k_2\alpha^2 & k_2\alpha^2 & k_3\alpha^2\\ k_4\alpha+k_5\alpha^2 & 1-k_6\alpha+k_5\alpha^2 & 0 \\ k_7+k_8\alpha+k_9\alpha^2 & k_{10}+k_9\alpha^2 & 1-k_1+k_{11}\alpha^2\end{smallmatrix}\end{bmatrix}, \\
\Upsilon_\alpha &\triangleq \begin{bmatrix}\color{black}\begin{smallmatrix}k_{12}\alpha^2\\ k_{13}\alpha^2 \\ k_{14}+k_{15}\alpha^2\end{smallmatrix}\end{bmatrix},
\end{align*}
where \textcolor{black}{$k_1\triangleq\frac{1-\bar{\sigma}^2}{2}$, $k_2\triangleq 24(n+4)n_c\varsigma\mathcal{L}^2$, $k_3\triangleq 2\varsigma$, $k_4 \triangleq \mathcal{L}^2/\chi$, $k_5\triangleq 12n(n+4)\mathcal{L}^2$, $k_6 \triangleq \chi$, $k_7\triangleq 12(n+4)n_s\varsigma\mathcal{L}^2(3+\bar{\sigma}^2)$, $k_8\triangleq 24(n+4)n_s\varsigma\mathcal{L}^2k_4$, $k_9\triangleq 24(n+4)n_s\varsigma\mathcal{L}^2(k_2+k_5)$, $k_{10}\triangleq 48(n+4)n_s\varsigma\mathcal{L}^2$, $k_{11}\triangleq \varsigma k_{10}$, $k_{12} \triangleq 24(n+4)n_c\varsigma G^2+6(n+4)^3n_c\varsigma\mathcal{L}^2\mu^2$, $k_{13} \triangleq 12n(n+4)G^2+3n(n+4)^3\mathcal{L}^2\mu^2$, $k_{14}\triangleq 48(n+4)n_s\varsigma G^2+12(n+4)^3n_s\varsigma\mathcal{L}^2\mu^2$ and $k_{15}\triangleq 24(n+4)n_s\varsigma\mathcal{L}^2(k_{12}+k_{13})$.}

Following the results in Lemmas~\ref{lemma:conensus}, \ref{lemma:optimality_gap} and \ref{lemma:grad_track_error}, and taking the total expectation on both sides, we obtain the following dynamical system:
\begin{align}
\Psi_{t+1} \leq \mathbf{M}_\alpha\Psi_t + \Upsilon_\alpha. \label{eq:dynamical_system}
\end{align}
Then the convergence of all players' actions to a neighborhood of the unique NE $\mathbf{x}^*$ of game $\Gamma$ is established based on the above dynamical system, as shown in the following theorem.

\begin{Theorem}\label{theorem:optimality}
Suppose Assumptions~\ref{assumption_graph}, \ref{assumption_local_f_lipschitz} and \ref{assumption_game_mapping} hold. Let \textcolor{black}{the auxiliary action variables $\{y^i_{jk,t}\}_{t\geq0}$,} the player's action $\{x^i_{j,t}\}_{t\geq0}$, and gradient tracker $\{\phi^i_{jk,t}\}_{t\geq0}$ be generated by \eqref{eq:algorithm} with a constant step-size $\alpha$ satisfying
\begin{align*}
\color{black}0<\alpha<\min\{\alpha_1,\alpha_2,\alpha_3,1/k_6,1\},
\end{align*}
where \textcolor{black}{$\alpha_1\triangleq[k_1^2k_6/(k_1k_2k_6+2k_1k_2k_4+4k_3k_4k_{10}+2k_3k_6k_7+2k_3k_6k_8)]^{\frac12}$, $\alpha_2\triangleq k_1k_4k_6/(k_1k_5k_6+2k_1k_4k_5)$, and $\alpha_3\triangleq \{k_1(2k_4k_{10}+k_6k_7+k_6k_8)/[k_1k_6k_9+2k_1k_4k_9+k_{11}(4k_4k_{10}+2k_6k_7+2k_6k_8)]\}^{\frac12}$.}
Then, $\rho(\mathbf{M}_\alpha)<1$, and we have 
\begin{itemize}
\item \textcolor{black}{$\sup_{\ell\geq t}\sum_{i=1}^N\sum_{k=1}^{n_i}\mathbb{E}[\|\mathbf{y}^i_{k,\ell}-\mathbf{1}_{n_i}\bar{y}^i_{k,\ell}\|^2]$} 
\item \textcolor{black}{$\sup_{\ell\geq t}\mathbb{E}[\|\bar{\mathbf{y}}_\ell - \mathbf{x}^*_\mu\|^2]$}
\item $\sup_{\ell\geq t}\sum_{i=1}^N\sum_{k=1}^{n_i}\mathbb{E}[\|\phi^i_{k,\ell} - \mathbf{1}_{n_i}\bar{\phi}^i_{k,\ell}\|^2]$
\end{itemize}
converge at a geometric rate with exponent $\rho(\mathbf{M}_\alpha)$. Moreover, 
\begin{align*}
&\color{black}\limsup_{t\to\infty}\sum_{i=1}^N\sum_{k=1}^{n_i}\mathbb{E}[\|\mathbf{y}^i_{k,t}-\mathbf{1}_{n_i}\bar{y}^i_{k,t}\|^2] \leq\mathcal{O}(\alpha^2),\\
&\color{black}\limsup_{t\to\infty}\mathbb{E}[\|\bar{\mathbf{y}}_t - \mathbf{x}^*_\mu\|^2]\leq\mathcal{O}(\alpha),
\end{align*}
\textcolor{black}{and we further have}
\begin{align*}
\color{black}\limsup_{t\to\infty}\mathbb{E}[\|\mathbf{x}_t - \mathbf{x}^*\|^2]\leq\mathcal{O}(\alpha)+\mathcal{O}(\mu).
\end{align*}
\end{Theorem}
\begin{Proof}
For the dynamical system \eqref{eq:dynamical_system}, if $\rho(\mathbf{M}_\alpha)<1$, then $\color{black}\mathbf{M}^t_\alpha$ converges to $\mathbf{0}$ at a geometric rate with exponent $\rho(\mathbf{M}_\alpha)$ \cite[Thm.~5.6.12]{Horn1990}, which implies that each element in $\Psi_t$ converges to some neighborhood of 0 with the same rate, respectively. 

\color{black}
The following lemma is leveraged to ensure $\rho(\mathbf{M}_\alpha)<1$:
\begin{Lemma}\label{lemma:matrix_spectral_radius}
(see \cite[Cor.~8.1.29]{Horn1990}) Let $A\in\mathbb{R}^{m\times m}$ be a matrix with non-negative entries and $\bm{\theta}\in\mathbb{R}^m$ be a vector with positive entries. If there exists a constant $\lambda\geq0$ such that $A\bm{\theta} < \lambda \bm{\theta}$, then $\rho(A) < \lambda$.
\end{Lemma}

To invoke Lemma~\ref{lemma:matrix_spectral_radius}, it suffices to set $\alpha<\frac1{k_6}$ to ensure each element of $\mathbf{M}_\alpha$ is non-negative.
Then, according to Lemma~\ref{lemma:matrix_spectral_radius}, it suffices to find a vector $\bm{\theta} \triangleq [\theta_1,\theta_2,\theta_3]^\top$ with $\theta_1,\theta_2,\theta_3>0$ such that $\mathbf{M}_\alpha\bm{\theta} < \bm{\theta}$, \textit{i.e.},
\begin{align*}
&(1-k_1+k_2\alpha^2)\theta_1+(k_2\alpha^2)\theta_2 + (k_3\alpha^2)\theta_3 < \theta_1,\\
&(k_4\alpha+k_5\alpha^2)\theta_1+ (1-k_6\alpha+k_5\alpha^2)\theta_2 < \theta_2,\\
&(k_7+k_8\alpha+k_9\alpha^2)\theta_1 + (k_{10}+k_9\alpha^2)\theta_2\\
&\quad\quad\quad\quad\quad\quad\quad\quad\quad\quad\quad\quad + (1-k_1+k_{11}\alpha^2)\theta_3 < \theta_3.
\end{align*}
Without loss of generality, letting $\theta_3 = 1$, it suffices to find $\theta_1$ and $\theta_2$ such that the following inequalities hold
\begin{subequations}\label{eq:ineq_t}
\begin{align}
&(k_2\theta_1+k_2\theta_2+k_3)\alpha^2 < k_1\theta_1, \label{eq:ineq_t1}\\
&(k_5\theta_1+k_5\theta_2)\alpha < k_6\theta_2-k_4\theta_1, \label{eq:ineq_t2}\\
&[k_9(\theta_1+\theta_2)+k_{11}]\alpha^2\nonumber\\
&\quad\quad\quad\quad\quad\quad<k_1-(k_7+k_8)\theta_1-k_{10}\theta_2,\label{eq:ineq_t3}
\end{align}
\end{subequations}
where we have forced $\alpha<1$ in the third inequality.
Letting $\theta_1=\frac{k_1k_6}{4k_4k_{10}+2k_6k_7+2k_6k_8}$ and $\theta_2 = \frac{k_1k_4}{2k_4k_{10}+k_6k_7+k_6k_8}$, we may solve \eqref{eq:ineq_t} that
\begin{align*}
&\alpha<\alpha_1,\alpha<\alpha_2,\alpha<\alpha_3,
\end{align*}
where $\alpha_1\triangleq[k_1^2k_6/(k_1k_2k_6+2k_1k_2k_4+4k_3k_4k_{10}+2k_3k_6k_7+2k_3k_6k_8)]^{\frac12}$, $\alpha_2\triangleq k_1k_4k_6/(k_1k_5k_6+2k_1k_4k_5)$, and $\alpha_3\triangleq \{k_1(2k_4k_{10}+k_6k_7+k_6k_8)/[k_1k_6k_9+2k_1k_4k_9+k_{11}(4k_4k_{10}+2k_6k_7+2k_6k_8)]\}^{\frac12}$,
which gives the desired range of $\alpha$.
\color{black}

To find out the steady-state value, taking the limsup on both sides of \eqref{eq:dynamical_system} gives
\begin{align*}
\limsup_{t\to\infty}\Psi_t \leq \mathbf{M}_\alpha\limsup_{t\to\infty}\Psi_t + \Upsilon.
\end{align*}
Notice that
\begin{align*}
\color{black}I_3 - \mathbf{M}_\alpha = \begin{bmatrix}\begin{smallmatrix}k_1-k_2\alpha^2 & -k_2\alpha^2 & -k_3\alpha^2\\ -k_4\alpha-k_5\alpha^2 & k_6\alpha-k_5\alpha^2 & 0 \\ -k_7-k_8\alpha-k_9\alpha^2 & -k_{10}-k_9\alpha^2 & k_1-k_{11}\alpha^2\end{smallmatrix}\end{bmatrix}.
\end{align*}
Its determinant can be obtained that
\begin{align*}
Det(I_3 - \mathbf{M}_\alpha) &= \color{black}\alpha(k_1-k_{11}\alpha^2)[k_1k_6-k_1k_5\alpha\\
&\quad\color{black}-k_2(k_4+k_6)\alpha^2].
\end{align*}
Then, we can obtain that
\begin{align*}
&\limsup_{t\to\infty}\sum_{i=1}^N\sum_{k=1}^{n_i}\mathbb{E}[\|\mathbf{y}^i_{k,t}-\mathbf{1}_{n_i}\bar{y}^i_{k,t}\|^2]\leq[(I_3 - \mathbf{M}_\alpha)^{-1}\Upsilon]_1\\
&\quad\color{black}\leq\frac{k_{12}\alpha^2\mathcal{O}(\alpha)+k_{13}\alpha^2\mathcal{O}(\alpha^2)+(k_{14}+k_{15}\alpha^2)\mathcal{O}(\alpha^3)}{\alpha(k_1-k_{11}\alpha^2)[k_1k_6-k_1k_5\alpha-k_2(k_4+k_6)\alpha^2]}\\
&\quad\color{black}=\mathcal{O}(\alpha^2),
\end{align*}
and
\begin{align*}
&\limsup_{t\to\infty}\mathbb{E}[\|\bar{\mathbf{y}}_t - \mathbf{x}^*_\mu\|^2]\leq[(I_3 - \mathbf{M}_\alpha)^{-1}\Upsilon]_2\\
&\quad\color{black}\leq\frac{k_{12}\alpha^2\mathcal{O}(\alpha)+k_{13}\alpha^2\mathcal{O}(1)+(k_{14}+k_{15}\alpha^2)\mathcal{O}(\alpha^3)}{\alpha(k_1-k_{11}\alpha^2)[k_1k_6-k_1k_5\alpha-k_2(k_4+k_6)\alpha^2]}\\
&\quad\color{black}=\mathcal{O}(\alpha).
\end{align*}
\color{black}Thus, it follows from \eqref{eq:x_minus_ybar} that
\begin{align*}
&\limsup_{t\to\infty}\mathbb{E}[\|\mathbf{x}_t - \mathbf{x}^*_\mu\|^2]\leq2\limsup_{t\to\infty}\mathbb{E}[\|\mathbf{x}_t - \bar{\mathbf{y}}_t\|^2]\\
&\quad\quad+2\limsup_{t\to\infty}\mathbb{E}[\|\bar{\mathbf{y}}_t - \mathbf{x}^*_\mu\|^2]\\
&\quad\leq2\limsup_{t\to\infty}\sum_{i=1}^N\sum_{k=1}^{n_i}\mathbb{E}[\|\mathbf{y}^i_{k,t}-\mathbf{1}_{n_i}\bar{y}^i_{k,t}\|^2]\\
&\quad\quad+2\limsup_{t\to\infty}\mathbb{E}[\|\bar{\mathbf{y}}_t - \mathbf{x}^*_\mu\|^2]=\mathcal{O}(\alpha).
\end{align*}
Applying Lemma~\ref{lemma:NE_gap}-(2), we have
\begin{align*}
&\limsup_{t\to\infty}\mathbb{E}[\|\mathbf{x}_t - \mathbf{x}^*\|^2]\leq2\limsup_{t\to\infty}\mathbb{E}[\|\mathbf{x}_t - \mathbf{x}^*_\mu\|^2]\\
&\quad\quad+2\|\mathbf{x}^*_\mu- \mathbf{x}^*\|^2=\mathcal{O}(\alpha)+\mathcal{O}(\mu),
\end{align*}
which completes the proof.\color{black}
\end{Proof}
\begin{Remark}
Theorem~\ref{theorem:optimality} shows that the players' actions converge to a neighborhood of the NE at a geometric rate with the error bounded by a term which is proportional to both the step-size \textcolor{black}{and the smoothing parameter}. When the step-size and the smoothing parameter are set smaller, the error bound decreases, leading to a better accuracy.
\end{Remark}
\begin{Remark}
The upper bound of the step-size depends on the global information, which may be infeasible to compute in a distributed mannar. However, the lower bound of the step-size is 0, which implies that the step-size can be set small enough to guarantee the convergence of the algorithm. The notion of a sufficiently small step-size is common in gradient tracking methods, see also \cite{Xi2016a,Xin2019b,Pu2020,Xin2019}.
\end{Remark}

\section{Numerical Simulations}\label{sec:simulation}

\begin{figure}[!t]
\centering
\includegraphics[width=2in]{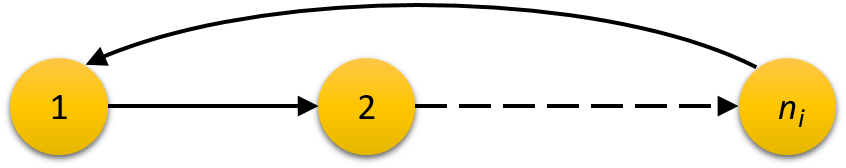}  %
\caption{Communication network.}
\label{fig:network.png}%
\end{figure}

\begin{figure}[!t]
\centering
\includegraphics[width=3.4in]{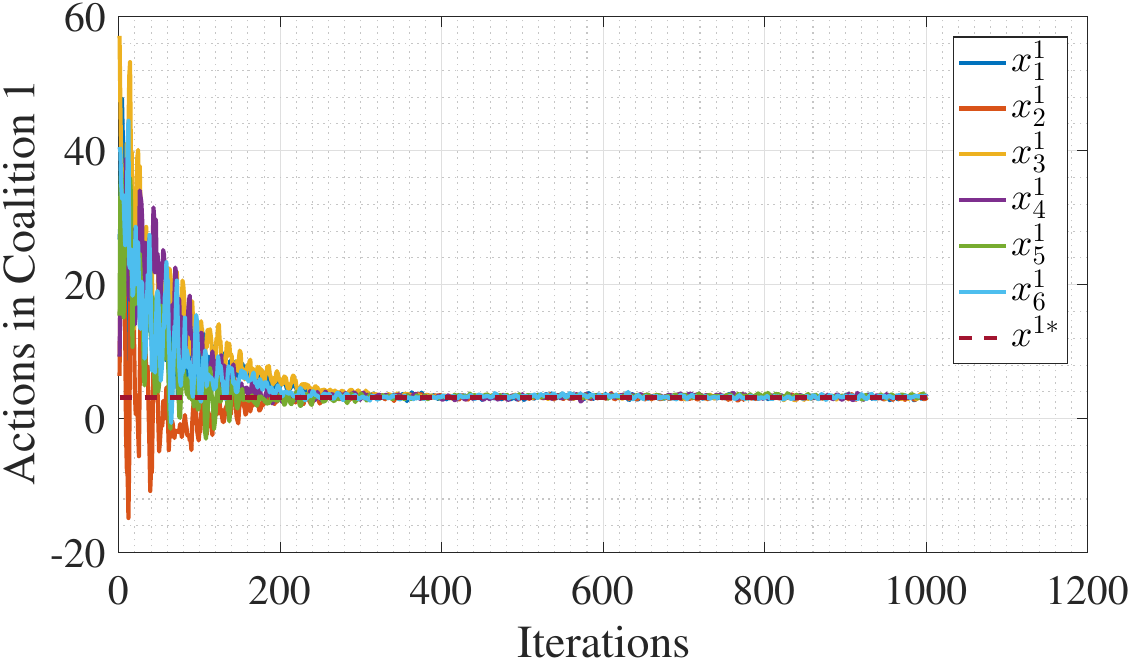}  %
\caption{Trajectories of players' actions in coalition 1.}
\label{fig:action_rgf_con1.png}%
\end{figure}

\begin{figure}[!t]
\centering
\includegraphics[width=3.4in]{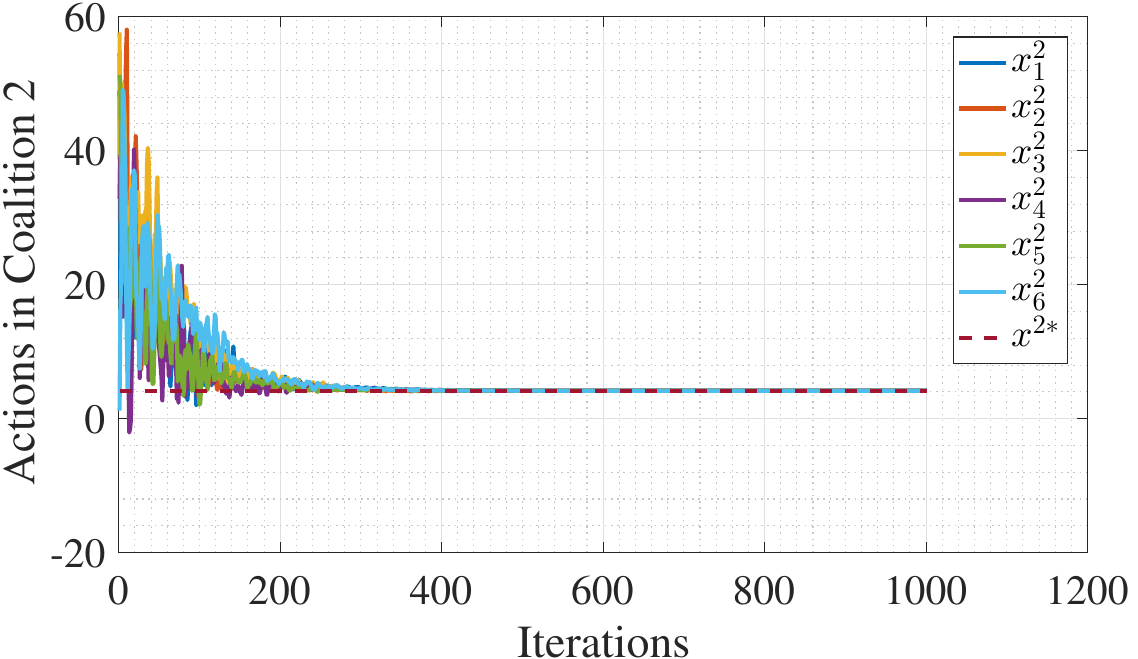}  %
\caption{Trajectories of players' actions in coalition 2.}
\label{fig:action_rgf_con2.png}%
\end{figure}

\begin{figure}[!t]
\centering
\includegraphics[width=3.4in]{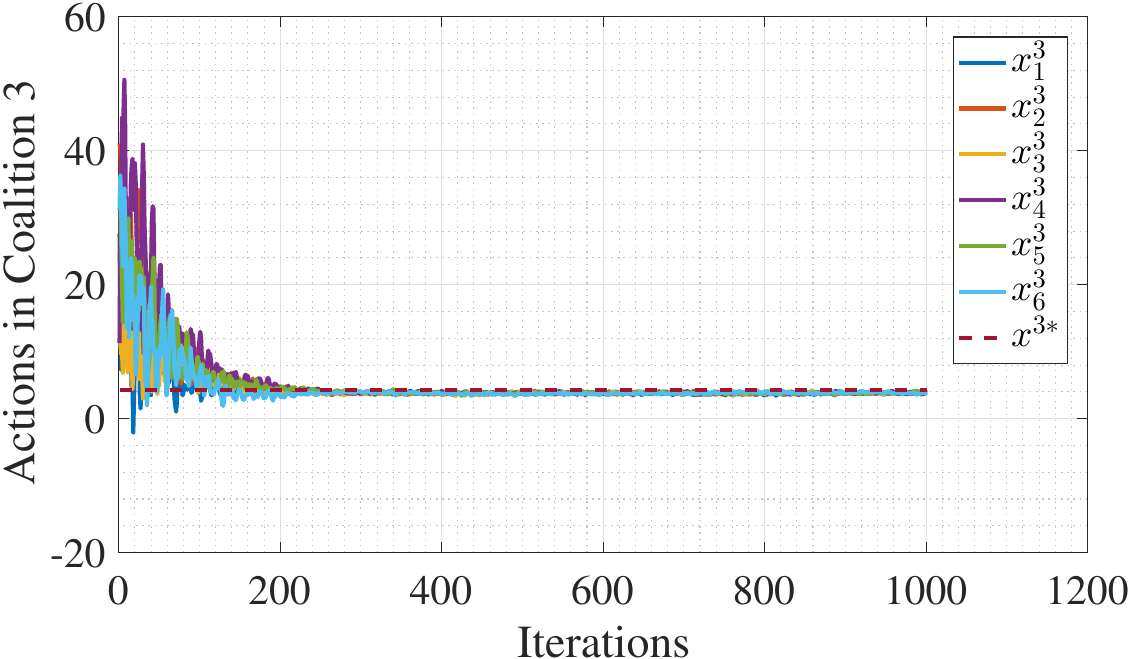}  %
\caption{Trajectories of players' actions in coalition 3.}
\label{fig:action_rgf_con3.png}%
\end{figure}

\begin{figure}[!t]
\centering
\includegraphics[width=3.4in]{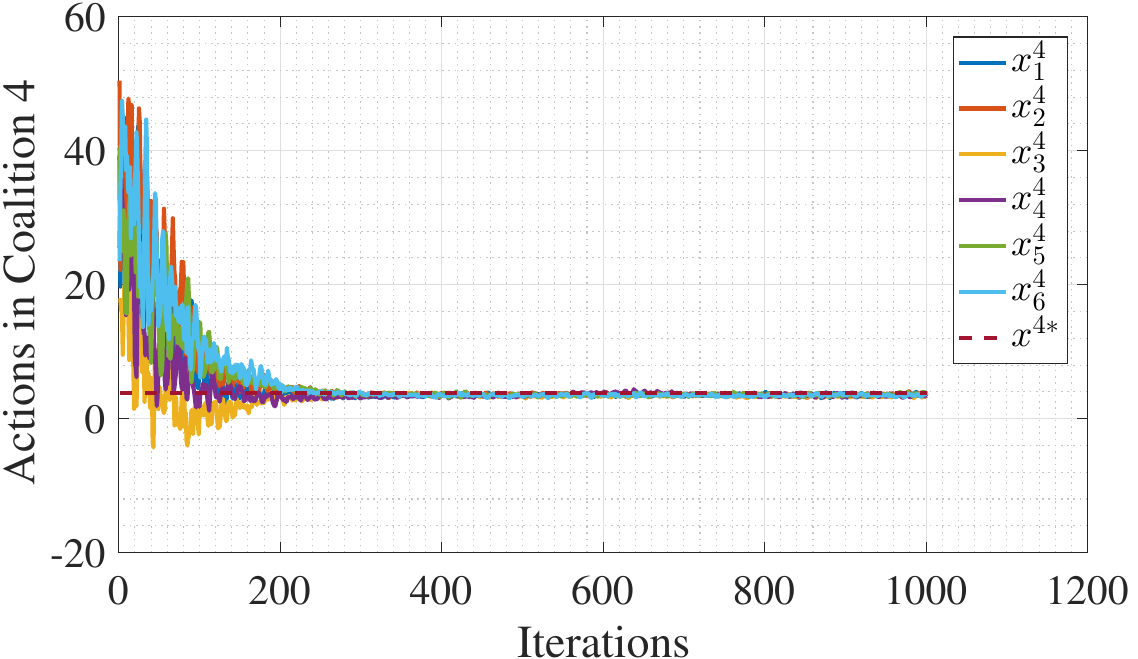}  %
\caption{Trajectories of players' actions in coalition 4.}
\label{fig:action_rgf_con4.png}%
\end{figure}

In this section, we illustrate the performance of the proposed algorithm with a numerical example for the Cournot competition game. The game consists of $N = 4$ coalitions with each coalition $i$ having $n_i = 6$ players. The local cost function of player $j \in \{1,\ldots,6\}$ in coalition $i\in\{1,\ldots,4\}$ is given by
$f^i_j(\mathbf{x}) = c^i_j(\mathbf{x}) - x^i_j{\color{black}r^i_j}(\mathbf{x})$,
where
$c^i_j(\mathbf{x}) = 5(x^i_j)^2 + 5x^i_j + \color{black}5(x^i_j - 6j)$, ${\color{black}r^1_j}(\mathbf{x}) = 60 - x^1_j-x^2_j-x^3_j-x^4_j$, ${\color{black}r^2_j}(\mathbf{x}) = 60 - x^2_j,{\color{black}r^3_j}(\mathbf{x}) = 60 - x^1_j-x^2_j$, and ${\color{black}r^4_j}(\mathbf{x}) = 60 - x^1_j-x^2_j-x^3_j$.
It is obvious to see that the problem is quadratic and hence satisfies Assumptions~\ref{assumption_local_f_lipschitz} and \ref{assumption_game_mapping}. The directed communication graph for each coalition $i$ is as shown in Fig.~\ref{fig:network.png}. Then, it is easy to find an associated adjacency matrix $A^i$ satisfying Assumption~\ref{assumption_graph}.
In the simulation, we set the smoothing parameter $\color{black}\mu = 10^{-4}$. To validate the convergence of the players' actions, we set the constant step-size $\alpha = 0.1$. The algorithm is initialized with arbitrary $x^i_{j,0}$, $y^i_{jk,0}$ and $\phi^i_{jk,0} = \pi^i_{jk}(\mathbf{x}_0)$. The trajectories of the players' actions for four coalitions are plotted in Figs.~\ref{fig:action_rgf_con1.png},\ref{fig:action_rgf_con2.png},\ref{fig:action_rgf_con3.png} and \ref{fig:action_rgf_con4.png}. As can be seen, all players' actions can approximately converge to the NE. 



\begin{figure}[!t]
\centering
\includegraphics[width=3in]{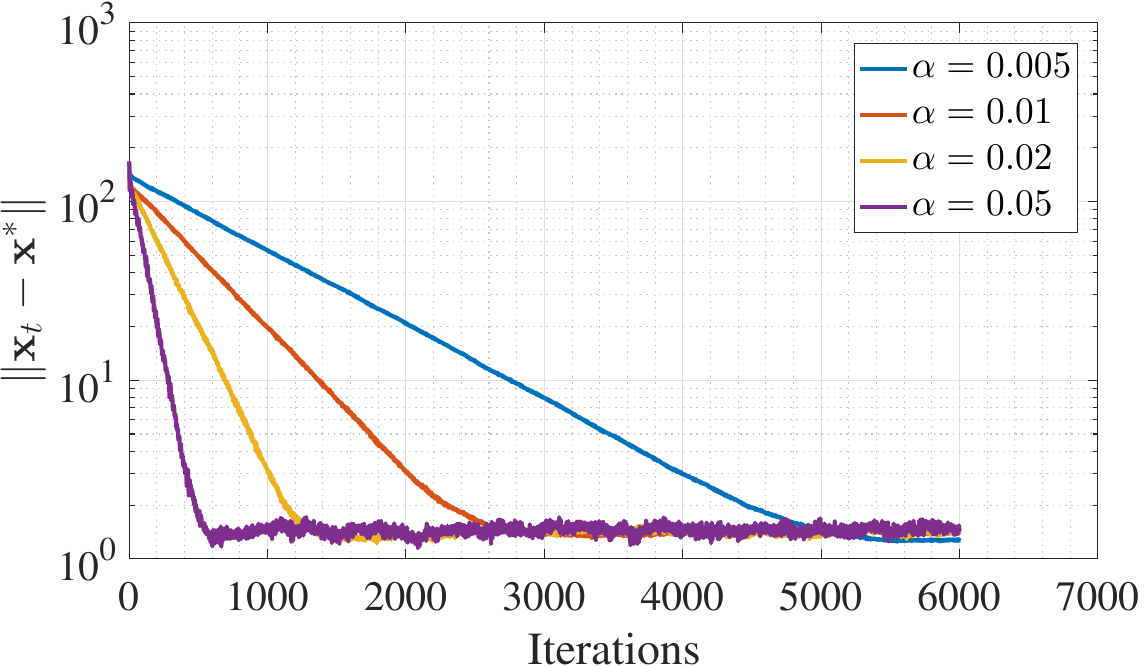}  %
\caption{Trajectories of the error gap $\|\mathbf{x}_t - \mathbf{x}^*\|$.}
\label{fig:convergence_rate.png}%
\end{figure}

To verify the convergence rate, we set the constant step-size $\alpha=0.005$, $0.01$, $0.02$ and $0.05$. The trajectories of the error gap $\|\mathbf{x}_t - \mathbf{x}^*\|$ with all these step-sizes are plotted in Fig.~\ref{fig:convergence_rate.png}. As can be observed, the error gap descends linearly in the log-scale plot for all cases. Moreover, when the step-size $\alpha$ is smaller, the convergence rate is slower but leading to a better accuracy, which verifies the derived results in Theorem~\ref{theorem:optimality}.

\section{Conclusions}\label{sec:conclusion}
In this paper, we have studied an $N$-coalition non-cooperative game problem, where players have no access to the explicit form but only the value of their local cost functions. A discrete-time gradient-free Nash equilibrium seeking algorithm, based on the gradient tracking method, has been proposed to search for the Nash equilibrium of the game. Under a strongly monotone game mapping condition, we have established that all players' actions converge linearly to a neighborhood of the Nash equilibrium with a sufficiently small constant step-size, where the gap is proportional to the step-size and the smoothing parameter. The performance of the proposed algorithm has been illustrated in numerical simulations.

\bibliographystyle{plain}
\bibliography{d_ne_coalition_rgf_reference}

\end{document}